\documentclass[11pt,a4paper]{article}

\usepackage[a4paper, left=1in, right=1in, top=1in, bottom=1in, includehead, includefoot, foot = 0pt, head=0pt]{geometry} 
\usepackage[group-separator={,},group-digits=integer]{siunitx} 
\usepackage{booktabs} 
\usepackage{enumitem} 

\usepackage[square,colon,numbers]{natbib}
\usepackage[pdftex,colorlinks,citecolor=blue,urlcolor=blue,linkcolor=red]{hyperref} 
\usepackage{graphicx} 
\parskip=\medskipamount 
\usepackage[small,bf,hang]{caption}	
\usepackage{amssymb} 
\usepackage{amsmath} 
\usepackage{bbm} 
\usepackage{amsbsy} 
\usepackage{mathrsfs} 
\usepackage[noend]{algorithm2e} 
\SetAlCapFnt{\small} 
\SetAlCapNameFnt{\small} 
\usepackage{multirow} 
\usepackage[dvipsnames,table]{xcolor} 
\usepackage{subcaption} 
\usepackage[hang,flushmargin]{footmisc} 
\usepackage{tikz} 
\usetikzlibrary{shapes,arrows,decorations.pathmorphing,backgrounds,positioning,fit,petri}
\usetikzlibrary{arrows.meta,
                chains,       
                positioning,
                quotes         
                }
\usetikzlibrary{shapes.multipart}                
\usetikzlibrary{matrix, fit} 
\usepackage{rotfloat}
\usepackage{titling} 
\makeatletter
\renewcommand{\paragraph}{%
	\@startsection{paragraph}{4}%
	{\z@}{0.5ex \@plus 1ex \@minus .2ex}{-1em}
	{\normalfont\normalsize\bfseries}%
}
\makeatother

\usepackage{authblk} 

\usepackage{fancyhdr}
\usepackage{fancybox}
\usepackage{tcolorbox}
\usepackage{latexsym,amsmath,xcolor,multicol,booktabs,calligra, adjustbox}
\usepackage{graphicx,listings,stackengine}
\usepackage{colortbl}
\definecolor{AXPBlue}{RGB}{0,52,102}
\definecolor{AXPGreen}{RGB}{152,203,2}
\usepackage{pstricks}
\usepackage{float}
\restylefloat{table}
\usepackage[hypcap=false]{caption}
\usepackage{adjustbox}
\usepackage{tikz}
\usepackage{mathtools}
\usepackage{mdframed}
\usepackage{amsthm}
\usepackage{scrextend}

\newtheoremstyle{break}
  {\topsep}{\topsep}%
  {\itshape}{}%
  {\bfseries}{}%
  {\newline}{}%
\theoremstyle{break}
 \usepackage{pstricks,pst-node,pst-tree}
 \usepackage{geometry}
 \usepackage{lscape}
 \definecolor{wheat}{rgb}{0.96, 0.87, 0.7}
 \definecolor{lightgray}{rgb}{0.9, 0.9, 0.9}
 \definecolor{aqua}{rgb}{0.0, 1.0, 1.0}
 \definecolor{mauve}{rgb}{0.88, 0.69, 1.0}
 \definecolor{skyblue}{RGB}{86,180,233}
 \definecolor{vermillion}{RGB}{213,94,0}
 \definecolor{purple}{RGB}{204,121,167}
 \definecolor{bluishgreen}{RGB}{0,158,115}

\newtheorem{definition}{Definition}
\newtheorem{remark}{Remark}

\usepackage{orcidlink}

\begin{document}
\sloppy

\title{\LARGE Including individual Customer Lifetime Value and competing risks in tree-based lapse management strategies}
\author[1,2]{Mathias Valla\thanks{Email: \href{mailto:mathias.valla@univ-lyon1.fr}{mathias.valla@univ-lyon1.fr} , URL: \url{https://mathias-valla.com}, ORCID: \orcidlink{0000-0003-4760-7849}\hspace{2mm} \href{https://orcid.org/0000-0003-4760-7849}{0000-0003-4760-7849}}} 
\author[3]{Xavier Milhaud\thanks{Email: \href{mailto:xavier.milhaud@univ-amu.fr}{xavier.milhaud@univ-amu.fr}, URL: \url{http://xaviermilhaud.fr}, ORCID: \orcidlink{0000-0002-3962-9434}\hspace{2mm} \href{https://orcid.org/0000-0002-3962-9434}{0000-0002-3962-9434}}}
\author[1]{Anani Olympio\thanks{Email: \href{mailto:anani.olympio@cnp.fr}{anani.olympio@cnp.fr}, ORCID: \orcidlink{0000-0003-1030-1125}\hspace{2mm} \href{https://orcid.org/0000-0003-1030-1125}{0000-0003-1030-1125}}}
\affil[1]{Univ Lyon, Université Claude Bernard Lyon 1, Institut de Science Financière et d'Assurances (ISFA), Laboratoire SAF EA2429, F-69366, Lyon, France.}
\affil[2]{Faculty of Economics and Business, KU Leuven, Belgium.}
\affil[3]{Aix-Marseille Univ, CNRS, Institut de Mathématiques de Marseille (I2M),UMR 7373, Marseille, France.}
\predate{}
\postdate{}
\date{}
\maketitle
\thispagestyle{empty}
\begin{abstract}
\noindent A retention strategy based on an enlightened lapse model is a powerful profitability lever for a life insurer. Some machine learning models are excellent at predicting lapse, but from the insurer’s perspective, predicting which policyholder is likely to lapse is
not enough to design a retention strategy. In our paper, we define a lapse management framework with an appropriate validation metric based on Customer Lifetime Value and profitability. \\
We include the risk of death in the study through competing risks considerations in parametric and tree-based models and show that further individualization of the existing approaches leads to increased performance. We show that survival tree-based models outperform parametric approaches and that the actuarial literature can significantly benefit from them. Then, we compare, on real data, how this framework leads to increased predicted gains for a life insurer and discuss the benefits of our model in terms of commercial and strategic decision-making. \\ [2mm]
\end{abstract}
\vfill


\section{Introduction}\label{sec:Introduction}
In life insurance, ``lapse risk" or ``persistency risk" is the risk that the policyholder will cancel the contract at a time other than when the issuer expected when pricing the contract (\cite{kpmg_2020}). A life insurance policy can lapse if the policyholder stops paying the premiums required to keep the policy in force. This can happen if the policyholder becomes unable or unwilling to make the premium payments or if the policyholder chooses to surrender the policy for its cash value. When a policy lapses, the coverage and benefits the policy provides are no longer in effect, and the policyholder will not receive any payout if they pass away after the policy has lapsed. This risk is not considered an insurance risk because the payment to the policyholder ``is not contingent on an uncertain future event that adversely affects the policyholder". However, lapse management is still undoubtedly a primary concern for life insurers. Lapses may substantially affect a company's solvency, its future profits and cash flows (\cite{Buchardt2014, Buchardt2015}) or its Asset and Liabilities Management (ALM) (\cite{kim2005,Gatzert2008,Eling2013,Eling2014}). The importance of measuring lapse and churn behaviours is global; it goes from yielding individual estimations of the Customer Lifetime Values (CLV) to being an estimator of a firm's profitability (\cite{Gupta2006b, GUPTA2009}) or strength (\cite{ascarza2018a}). Therefore, this paper focuses on developing strategies to prevent lapses before they occur: for a life insurer, an enlightened and proactive lapse management strategy (LMS) is critical for successful monitoring and steering. This paper is about defining a framework for a life insurer to measure and optimise the future loss or profit to be expected when applying such an LMS.\\
\\
Part of the literature on lapse management adopts an economic-centred point of view (\cite{Dar1989, Kuo2003, Kagraoka2005, Cox2006, Kiesenbauer2012, Russel2013, sirak2015, vasudev2016a, Nolte2017, Poufinas2018, Yu2019}); we refer the reader to the complete bibliometric analysis on this topic by \cite{Shamsuddin2022} for a summarised view of all these references. This economic-centred research aims to determine lapse factors like interest rates, gross domestic product, or unemployment rates. They are driven by economic hypotheses such as the emergency fund hypothesis (lapsing is a way of constituting an emergency fund), the policy replacement hypothesis (lapsing will occur when one changes its policy) or the interest rate hypothesis (lapsing depends strongly on rate change and arbitration). \\
On the other hand, a large part of the literature investigates the individual determinants of lapse with policyholder-centred approaches. Micro-oriented features such as policyholder's personal information or the policy characteristics have shown to give valuable insights into lapse behavior (\cite{Renshaw1986,milhaud2011,Eling2014,Hwang2021}). \cite{Curak2015} as well as \cite{Gemmo2016}'s works indicate that policyholders' features such as age and the number of beneficiaries are significant lapse factors, whereas \cite{sirak2015} dismissed those results. A recent work from \cite{loisel_piette_tsai_2021} proposes a comparison of lapse management strategies based on an innovative evaluation metric derived from the Customer Lifetime Value (CLV). \cite{ hu_ohagan_sweeney_ghahramani_2021} investigates the benefits of incorporating spatial analysis in lapse modeling, and \cite{AZZONE2022116261} shows with an approach based on random forests that microeconomic features such as the company's commercial approach for instance -  is determining in the lapse decision. In contrast, macro-economical features only have a limited effect. This variety of results – sometimes contradicting each other – demonstrates the active interest in this research problem.\\
\\
This paper focuses on lapse management strategy and retention targeting and will contribute to the existing literature on the relationship between retention strategy and lapse prediction: as in \cite{ascarza2018a} and \cite{loisel_piette_tsai_2021}, our goal is not only to model the lapse behavior but rather to select policyholders that are expected to generate future profit, if targeted by a retention strategy. This work shows that a well-chosen strategy, based on individualized CLV and directed towards a well-chosen target, increases the insurer's expected profitability. A critical concept that motivates many CLV-driven decisions is that customers should be judged as assets based on their future profitability for the insurer. Thus, since retention often serves as the basis for CLV models (\cite{Gupta2006a,Donkers2007,Lemmens2020} - sometimes specifically designed for targeting tasks (\cite{VONMUTIUS2021}) - and since CLV considerations should drive retention management, it seems natural to extend the existing life insurance applications linking those topics together. We make decision-making a central concern of our work and suggest proactive lapse management tools allowing the insurer to undertake actions to prevent the causes of lapse; that is opposed to a reactive management approach where decisions are taken after lapses and aim at recapturing lost policyholders.\\
\\
The goal of this paper is to create an individualized CLV model that will be used to enhance classical binary churn models. We will then have a model for lapse management strategy and retention targeting that we further improve with tree-based survival analysis and competing risks considerations. The global framework is directly inspired by \cite{loisel_piette_tsai_2021}. We try in this paper to build from that existing work and extend it. We model an individual future CLV with a new survival approach for which the risks of death and lapse are treated as mutually exclusive competing risks. For this purpose, we introduce parametric approaches - Cox cause-specific and subdistribution models - as well as tree-based survival models - Random Survival Forest (RSF) and Gradient boosting survival analysis. We focus here on tree-based models as they are often considered state-of-the-art models (\cite{Grinsztajn2022}). Thus we introduce tree-based machine learning algorithms for binary prediction, including Classification and Regression Tree (CART), Random forests (RF), and Extreme Gradient Boosting (XGBoost) to lapse behavior modeling. CART and XGBoost (\cite{milhaud2011,loisel_piette_tsai_2021}) were used in the literature for lapse modeling but have yet to be applied to predicting life insurance lapses in a competing risk setting. To our knowledge, while Random Survival Forest has been used for churn prediction recently (\cite{Routh2021}), both RSF and Gradient boosting survival analysis have never been used for that purpose before in an actuarial context. \\
Our contribution to the actuarial literature is twofold. First, we detail a two-step lapse management modeling approach:  we fit parametric and tree-based competing risk individual survival models to estimate individualized future CLVs that are part of an evaluation metric for tree-based lapse management models. Second, this work includes a business-oriented discussion of the results achieved by this framework, which is missing from existing similar approaches.\\
The results and discussions show that a CLV-based lapse management strategy very often outperforms a more classical binary classification approach, even with competing risks and individualized considerations. When the latter yields profitable retention gain, the former can produce higher profits, up to more than 60\%. If a loss-inducing retention strategy is considered, our methodology limits the loss considerably, often setting 0 as a floor limit or even turning it into a profit-inducing retention strategy. Sensitivity analysis explores the influence of conjectural and structural parameters.\\
\\
The rest of this paper is structured as follows. We briefly outline the data used in our study in Section~\ref{sec:Data}. In Section~\ref{sec:Framework}, we then introduce the binary classification models we selected and detail our study's methodology, describing the classical and CLV-based performance measures and discussing substantial parametrization improvements over existing approaches. Then, Section~\ref{sec:Methodology} details our two-step methodology, with the parametric and non-parametric modelings of individual survival predictions, in a competing risks framework and then their implementation in the tree-based classification approaches considered. Section~\ref{results} presents the real-life application we considered and the different results it produces. Those results are analyzed and discussed in Section~\ref{Discussion} with commercial and strategical decision-making orientations. Eventually, Section~\ref{Conclusion} concludes this paper.

\section{Data}\label{sec:Data}
We apply our framework to a real-world insurance portfolio. For privacy reasons, all the data, statistics, product names and perimeters presented in this paper have been either anonymized or modified. All analyses, discussions and conclusions remain unchanged.\\
\\
We illustrate our methodology with a life insurance portfolio
from a French insurer contracted between 1997 and 2018. Each record in the data set represents a unique policy for a unique policyholder. In the following sections, we will often refer to a unique pair of policy and policyholder by the term ``subject". The dataset contains 251,325 rows with 248,737 unique policies and 235,076 unique policyholders. It means that some policies are shared between several policyholders and that one individual can detain several insurance policies. The dataset contains 43 covariates described in Table~\ref{tab:dataset_description}.\\

\begin{table}[H]
  \centering
  \caption{Data set description}
  \begin{adjustbox}{width=\textwidth}
    \begin{tabular}{lll}
    \toprule
    \multicolumn{2}{c}{Covariates (\textcolor{teal}{Numerical}, \textcolor{purple}{Categorical}, \textcolor{orange}{Date})} & \multicolumn{1}{c}{Description} \\
    \midrule
    \multicolumn{1}{c}{\multirow{2}[2]{*}{ID}} & \textcolor{purple}{CDI\_ID\_PERSONNE} & Policyholder (PH) unique ID \\
          & \textcolor{purple}{CDI\_ID\_CONTRAT} & Policy unique ID \\
\cmidrule{1-1}    \multicolumn{1}{c}{\multirow{9}[2]{*}{PH-level information}} & \textcolor{orange}{CDI\_DT\_NAISSANCE} & PH's birth date (main PH when several policyholders owns one policy) \\
          & \textcolor{teal}{\underline{Age\_souscription}} & PH's age at subscription \\
          & \textcolor{teal}{\underline{Nb\_Contrats}} & Number of different policies owned by the policyholder \\
          & \textcolor{purple}{\underline{CDI\_CD\_SEXE}} & PH's gender (1=Female; 2=Male; other=Non precised) \\
          & \textcolor{purple}{CDI\_DESTINATAIRE\_COURRIER} & Anonymised PH's name \\
          & \textcolor{purple}{CDI\_NUM\_ET\_NOM\_VOIE} & Anonymised PH's address \\
          & \textcolor{teal}{CDI\_CD\_POSTAL} & Anonymised PH's postcode \\
          & \textcolor{purple}{CDI\_COMMUNE} & Anonymised PH's city of residence \\
          & \textcolor{purple}{CDI\_TOP\_ASSURE} & Binary: 1 if PH is the main PH on the policy, 0 otherwise \\
\cmidrule{1-1}    \multicolumn{1}{c}{\multirow{10}[2]{*}{Policy-level informations}} & \textcolor{purple}{CDI\_TYPE\_PRODUIT} & Type of product (``Top-end product" or ``Classical product") \\
          & \textcolor{purple}{\underline{CDI\_NOM\_PRODUIT}} & Name of life insurance product (``Product 1", ``Product 2" or ``Product 3")\\
          & \textcolor{purple}{CDI\_PARTENAIRE} & Name of the insurance distributor \\
          & \textcolor{orange}{CDI\_DATE\_DEB\_CONTRAT} & Policy's start date \\
          & \textcolor{orange}{CDI\_DATE\_FIN\_CONTRAT} & Policy's end date \\
          & \textcolor{teal}{\underline{START\_YEAR}} & Policy's start year \\
          & \textcolor{teal}{END\_YEAR} & Policy's end year \\
          & \textcolor{teal}{\underline{SENIORITY}} & Policy's seniority (final seniority if the policy is ended, current seniority otherwise) \\
          & \textcolor{purple}{STATE} & Policy's state (``Active", ``Lapsed", or ``Death" if the policy ended following PH's death) \\
          & \textcolor{teal}{YEAR}  & Last year of observation \\
\cmidrule{1-1}    \multicolumn{1}{c}{External data} & \textcolor{teal}{DISCOUNT RATE} & Discount rate corresponding to the last year of observation \\
\cmidrule{1-1}    \multicolumn{1}{c}{\multirow{15}[2]{*}{Policy's cumulated financial flows}} & \textcolor{teal}{\underline{TOTAL\_PREMIUM\_AMOUNT}} & Total face amount of the policy \\
          & \textcolor{teal}{TOTAL\_EURO\_PREMIUM\_AMOUNT} & Face amount of the policy in euros \\
          & \textcolor{teal}{TOTAL\_UC\_PREMIUM\_AMOUNT} & Face amount of the policy in units of account \\
          & \textcolor{teal}{ARBITRATION\_EURO} & Cumulated arbitration amount of the policy in euros \\
          & \textcolor{teal}{ARBITRATION\_UC} & Cumulated arbitration amount of the policy in units of account \\
          & \textcolor{teal}{FEES\_EURO} & Cumulated fees amount of the policy in euros \\
          & \textcolor{teal}{FEES\_UC} & Cumulated fees amount of the policy in units of account \\
          & \textcolor{teal}{OTHER\_EURO} & Cumulated other parts of the face amount of the policy in euros \\
          & \textcolor{teal}{OTHER\_UC} & Cumulated other parts of the face amount of the policy in units of account \\
          & \textcolor{teal}{PREMIUM\_EURO} & Cumulated payments amount of the policy in euros \\
          & \textcolor{teal}{PREMIUM\_UC} & Cumulated payments amount of the policy in units of account \\
          & \textcolor{teal}{PROFIT SHARING\_EURO} & Cumulated profit sharing amount of the policy in euros \\
          & \textcolor{teal}{PROFIT SHARING\_UC} & Cumulated profit sharing amount of the policy in units of account \\
          & \textcolor{teal}{CLAIM\_EURO} & Cumulated partial or total lapsed amount of the policy in euros \\
          & \textcolor{teal}{CLAIM\_UC} & Cumulated partial or total lapsed amount of the policy in units of account \\
\cmidrule{1-1}    \multicolumn{1}{c}{\multirow{5}[2]{*}{Covariates derived from financial flows}} & \textcolor{teal}{\%TOTAL\_UC\_PREMIUM\_AMOUNT} & Percentage of the face amount in units of  account \\
          & \textcolor{teal}{\%TOTAL\_EURO\_PREMIUM\_AMOUNT}  & Percentage of the face amount in euros \\
          & \textcolor{teal}{\%CLAIM\_UC} & Percentage of the  face amount in units of account that was lapsed \\
          & \textcolor{teal}{\%CLAIM\_EURO} & Percentage of the  face amount in euros that was lapsed \\
          & \textcolor{teal}{\%CLAIM}  & Percentage of the total face amount that was lapsed \\
\cmidrule{1-1}    Target covariate & \textcolor{purple}{\underline{EVENT}} & Policy's state (0=Active, 1=Lapsed, 2 ended following PH's death) \\
    \end{tabular}%
    \end{adjustbox}
  \label{tab:dataset_description}%
\end{table}%

\begin{figure}[!htb]
    \centering
    \begin{minipage}[t]{0.3\linewidth}
        \centering
        \includegraphics[width=1\linewidth]{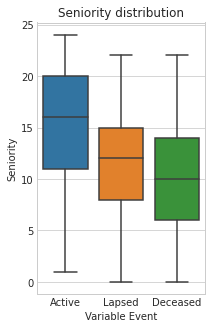}
    \end{minipage}
    \begin{minipage}[t]{0.65\linewidth} 
        \centering
        \includegraphics[width=1\linewidth]{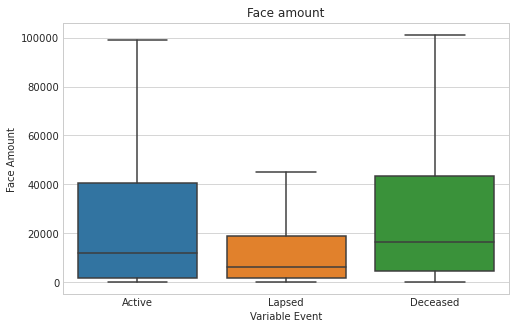}
    \end{minipage}
    \caption{Seniorities and face amounts distributions}
    \label{fig:seniority_distribution}
\end{figure}

The data set represents policies that are majority owned by men (57.4\%) for a mean censored seniority time of 13.4 years. Three products are present in the dataset. Product one was chosen by 72\% of policyholders, product 2 by 25\% and product 3 by 3\%.\\
Regarding their state, 61\% of the policies are still active, 22\% lapsed, and 17\% ended after the PH's death. We chose here to present the distribution of the variable \emph{SENIORITY} as it is the response variable in our survival models. Its modeling has a critical influence on CLV, thus, on our lapse management strategy framework. We also chose to show the distribution of the variable \emph{TOTAL PREMIUM AMOUNT} representing the most recent observed face amount for every subject, as it is a known determinant of lapse behavior. We are aware that this covariate is a rather dynamic one as its value is updated at every payment, total or partial lapse, profit sharing, arbitration or even fees movements on a policy, and only considering its most recent value ignores a large part of the insights it can provide. Without any better option, we can only use \emph{TOTAL PREMIUM AMOUNT} as it is and defer any dynamic considerations for future work.\\
\\
The seniorities and most recent face amount recorded before the potential end of the policy are distributed as in Figure~\ref{fig:seniority_distribution}:\\
\newline
Without further analyzing the data, we can note several things. First, we can see that the mean censored seniority of 13.4 years is not equally distributed among our subjects. Active contracts tend to be older than lapsed ones, themselves older than policies that ended with the policyholder's death. That emphasizes the importance of several contributions, and the apparent difference in seniority regarding the cause of the policy's termination encourages a competing risks approach to analyze survival. Moreover, if we suspect lapse and death to be highly dependent on individual characteristics - such as the policyholder's age - this also supports an individualized survival analysis. Eventually, we can see that the last face amount observed is significantly lower for lapsed policies. It confirms our first intuitions and the face amount will be included in our model.\\
\\
Among the covariates introduced in Table~\ref{tab:dataset_description}, several play a central part in our two-step modeling approach. First, the competing risks survival analysis step where \emph{SENIORITY} will be the response variable, and all other covariates, including individual data and financial flows, are potential explanatory variables. The binary classification second step aims at predicting the \emph{EVENT} outcome with minor transformations explained in Section~\ref{sec:Framework} below. It is equivalent to using \emph{STATE} as a target variable, as they are entirely similar. As a result, all covariates are not utilized and our predictions are solely based on the covariates \underline{underlined} in Table~\ref{tab:dataset_description} as they appeared to be of interest to insurers.


\section{Framework}\label{sec:Framework}
This section describes a modeling approach that follows \cite{loisel_piette_tsai_2021}'s work. Our contributions place our work in a framework that differs from it by being only future-oriented, by a precise and individualized analysis of retention probabilities and by choosing a classification framework instead of regression. We chose to use a majority of their existing notations here.\\
Usual lapse management models based on classification aim to predict whether a policyholder will lapse. They may perform very well at that specific task, but it only reflects some aspects of this economic problem. Indeed, the literature is clear (\cite{ascarza2018a}), and many policyholders may be predicted as ``lapsers" but may not be profitable to the insurance company if targeted. In that case, keeping such policyholders would be irrational, and an efficient model should not predict them as targets. Targeting policyholders is an economic problem that requires an economic measure to assess. We propose to consider a measure based on the discounted expected profit of all the policies, in other words, the sum of all $(\mathrm{CLVs})$. Optimizing a lapse, churn, or other prediction tasks with business-related measures is not new. However, to our knowledge, none of the existing approaches uses individualized future CLVs and models the profit of retention strategy by accounting for competing risks or using survival tree-based models.\\
\\
$\mathrm{CLV}$ is a well-studied subject in marketing and business economics and has also been modeled in an insurance context. For a given subject $i$, her future $\mathrm{CLV}$ at horizon $T$ can be modeled as 

\begin{equation}
    \prescript{F}{}CLV_{i}\left(\boldsymbol{p}_{i}, \boldsymbol{F}_{i}, \boldsymbol{r}_{i}, \boldsymbol{d}, T\right)=\sum_{t=0}^{T} \frac{p_{i, t} F_{i, t} r_{i, t}}{\left(1+d_{i, t}\right)^{t}} ,
\end{equation}

with $t$ in years, $t=0$ represents the last observation point for subject $i$. The quantity $p_{i, t}$ is her profitability ratio as a proportion of $F_{i, t}$, representing her face amount observed at time $t$. The quantity $r_{i, t}$ is the $i$-th subject's probability of still being active at time $t$, and naturally, $d_{i, t}$ is the discount rate at time $t$, for subject $i$. We argue that both the profitability ratio and the discount rate should be as individualized as possible - either at the product or policy level - as $\prescript{F}{}CLV_{i}$ reflects the individualized risk of policyholder $i$ to the insurer. It is also worth mentioning that evaluating discount rates is well beyond the scope of this paper as it is complex and subject to significant judgment; for further details, we refer the astute reader to a variety of papers on the subject (\cite{Burrows2000RiskDR}, \cite{Oh2018}, \cite{blum2019discount}).\\
\\
It is worth pointing out that $\prescript{F}{}CLV_{i}$ does not represent the global profit generated by subject $i$ from her policy's first year until time $T$ as in \cite{loisel_piette_tsai_2021}; it rather represents the future $T$ years of profit. $\prescript{F}{}CLV_{i}$ is not to be compared with the Cash Surrender Value but rather with the Fair Market Value (FMV) of the outstanding liabilities. The only difference with the latter is that $\prescript{F}{}CLV_{i}$ is based on the insurer's knowledge of its portfolio, thus computed with its own profitability and discount parameters rather than with market-consistent considerations. Whereas the portfolio market value would consist of financial instruments that replicate the insurance liability cash flows. In our framework, the life insurer is more interested in maximizing its own realistic profitability rather than a sum of individual market values.\\
We suggest a model for the insurer's estimated profit - or loss - resulting from a lapse management strategy (LMS). In order to do that, we will compare the expected value of the portfolio before and after applying a given strategy. We are aware that there could be infinite ways to design a retention campaign: offering a punctual incentive, recurrent services or more profit sharing, for instance. Here, we define what we will consider an LMS.
\vspace{3 mm}
\begin{mdframed}[
  leftmargin=0.5\parindent,
  rightmargin=0.5\parindent,
  skipabove=\topsep,
  skipbelow=\topsep,
  hidealllines=true
  ]
\begin{definition}[Lapse management strategy]\label{LMS} A T-years lapse management strategy is modeled by the offer of an incentive $\delta_{i}$ to subject $i$ if she is targeted. The incentive is expressed as a percentage of her face amount and should not exceed the profitability ratio $p_{i,t}$, at any time point $t$. Contacting the targeted policyholder has a fixed cost $c$ and after contact, the incentive is accepted with probability $\gamma_{i}$.  A targeted subject who accepts the incentive will be considered as an ``acceptant" who will never lapse. In our dataset, any subject that has never been observed to lapse is considered as an ``acceptant" and her probability of being active at year $t \in [0, T]$ is denoted $r_{i,t}^{\text {acceptant}}$.  Conversely, a subject who refuses the incentive and prefers to lapse will be considered as a ``lapser". In our dataset, a subject is labeled as a lapser whenever she has been observed to lapse at year $t=0$, and her probability of being active at year $t$ is denoted $r_{i,t}^{\text {lapser}}$. A lapse management strategy is uniquely defined by the parameters $(\boldsymbol{p},\boldsymbol{\delta},\boldsymbol{\gamma}, c, T)$
\end{definition}
\end{mdframed}
\vspace{3 mm}
It is to be noted that even if the framework involves a time dimension, it is still a static approach: the insurer would run all analyses on its portfolio at one given time and apply an appropriate LMS immediately.\\
Even if this definition is already a simplification of any real-life insurance retention strategy, various constraints and the data and tools at the insurer's disposal do not always allow to conduct such a study. In the following section, we consider a simplified version of this framework by assuming that $p_{i, t}, F_{i, t}$, and $d_{i,t}$ remain constant across time, and denoted $p_{i}$, $F_{i}$ and $d_{i}$ hereafter, with $F_{i}$ being the most recent face amount observed for subject $i$. Moreover, we set $\gamma_{i}$ and $\delta_{i}$ to be the same for all subjects and denoted as $\gamma$ and $\delta$ hereafter. The constraint that $\delta < min(p_{i})$ detailed in Definition~\ref{LMS} still holds. Finally, the last observed state of subject $i$ is denoted $y_i$, with $y_i=1$ if the policy is lapsed, $y_i=0$ otherwise. \\
\\
With those considerations, we can then define the control portfolio's future value as  
\begin{equation}
    \begin{aligned}
    \prescript{F}{}CPV&(\boldsymbol{p},\delta,\gamma, c, T)=\\
    &\phantom{aa} \displaystyle \sum_{i=1}^{n} \prescript{F}{}CLV_i\left({p}_{i}, F_i, \boldsymbol{r}^{\text {acceptant}}_{i}, d_{i}, T\right) \cdot \mathbf{1}\left(y_{i}=0\right) \\
    &+ \displaystyle \sum_{i=1}^{n} \prescript{F}{}CLV_i\left({p}_{i}, F_i, \boldsymbol{r}^{\text {lapser}}_{i}, d_{i}, T\right) \cdot \mathbf{1}\left(y_{i}=1\right) .
    \end{aligned}
\end{equation}
It represents the hypothetical value of the portfolio, considering that: 
\begin{itemize}
  \item every subject that did not lapse up to her last observation point - $y_i = 0$ at $t=0$ - has a vector of retention probabilities of $\boldsymbol{r}^{\text {acceptant}}_{i}$;
  \item every subject that has been observed to lapse - $y_i = 1$ at $t=0$ - has a vector of retention probabilities of $\boldsymbol{r}^{\text {lapser}}_{i}$
\end{itemize}

\begin{mdframed}[hidealllines=true]
\begin{remark}\label{remark_1}
It is important to note that this does not reflect the actual future value of the portfolio - as the future CLV of lapsers should be $0$ - but rather its hypothetical expected future value given the nature (lapser or not) of every subject but not their actual states (actually lapsed or not). It represents this hypothetical future CLV of all subjects if no customer relationship management about lapses is carried out.
\end{remark}
\end{mdframed}

A classification algorithm would take the lapse indicator $y_i$ as a target variable and yield predictions $\hat{y}_{i}$. Given a lapse management strategy and such a classification algorithm, we define the lapse managed portfolio future value by\\

\begin{equation}
\scalebox{0.9}{\parbox{1\linewidth}{%
$$
\begin{aligned}
\prescript{F}{}L&MPV(\boldsymbol{p},\delta, \gamma, c, T)=\\
& \phantom{aaaaaaaiii} \displaystyle \sum_{i=1}^{n} \prescript{F}{}CLV_i\left({p}_{i}, F_i, \boldsymbol{r}^{\text {acceptant}}_{i}, d_{i}, T\right)  \cdot \mathbf{1}\left(y_{i}=0, \hat{y}_{i}=0\right)\\
& \phantom{aaaaaiii} + \displaystyle \sum_{i=1}^{n} \prescript{F}{}CLV_i\left({p}_{i}, F_i, \boldsymbol{r}^{\text {lapser}}_{i}, d_{i}, T\right)  \cdot \mathbf{1}\left(y_{i}=1, \hat{y}_{i}=0\right)\\
& \phantom{aaaaaiii} + \displaystyle \sum_{i=1}^{n} \prescript{F}{}CLV_i\left({p}_{i}-\delta, F_i, \boldsymbol{r}^{\text {acceptant}}_{i}, d_{i}, T\right)  \cdot \mathbf{1}\left(y_{i}=0, \hat{y}_{i}=1\right)\\
& \phantom{aaaai} + \gamma \cdot \displaystyle \sum_{i=1}^{n} \prescript{F}{}CLV_i\left({p}_{i}-\delta, F_i, \boldsymbol{r}^{\text {acceptant}}_{i}, d_{i}, T\right)  \cdot \mathbf{1}\left(y_{i}=1, \hat{y}_{i}=1\right)\\
& + (1-\gamma) \cdot \displaystyle \sum_{i=1}^{n} \prescript{F}{}CLV_i\left({p}_{i}, F_i, \boldsymbol{r}^{\text {lapser}}_{i}, d_{i}, T\right)  \cdot \mathbf{1}\left(y_{i}=1, \hat{y}_{i}=1\right)\\
& \phantom{aaaaaiii} - \displaystyle \sum_{i=1}^{n} c \cdot \mathbf{1} \left(\hat{y}_{i}=1\right) .
\end{aligned}
$$
}}
\end{equation}

Clearly, the sums appearing in the formulas above could be grouped to make them more concise. We chose not to do so for the sake of visualization: we can distinctly see each possible case in each summand.\\
Then, we define the economic metric of the algorithm as the retention gain, the future profit generated by the retention management strategy over $T$ years as\\
\begin{equation}
    RG(\boldsymbol{p}, \delta, \gamma, c, T)=\operatorname{\prescript{F}{}LMPV}(\boldsymbol{p}, \delta, \gamma, c, T)-\prescript{F}{}CPV(\boldsymbol{p}, \delta, \gamma, c, T) ,
\end{equation}

which can be simplified as\\
\begin{equation}
\begin{aligned}
RG(\boldsymbol{p}, \delta, \gamma, c, T)=&\displaystyle \sum_{i=1}^{n} \Biggr[ \gamma \Big[ \prescript{F}{}CLV_i\left({p}_{i}-\delta, F_i, \boldsymbol{r}^{\text {acceptant }}_i, d_{i}, T\right)\\
&-\prescript{F}{}CLV_i\left({p}_{i}, F_i, \boldsymbol{r}^{\text {lapser }}_i, d_{i}, T\right)\Big] \cdot \mathbf{1}\left(y_{i}=1, \hat{y}_{i}=1\right)\\
&-\prescript{F}{}CLV_i\left(\delta, F_i, \boldsymbol{r}^{\text {acceptant }}_i, d_{i}, T\right) \cdot \mathbf{1}\left(y_{i}=0, \hat{y}_{i}=1\right)\Biggr]\\
& - \displaystyle \sum_{i=1}^{n} c \cdot \mathbf{1} \left(\hat{y}_{i}=1\right) .
\end{aligned}
\end{equation}

This evaluation metric can now be derived into an individual retention gain measure. More specifically, we define $z_{i}$ as\\
\begin{equation}
  z_{i} =
    \begin{cases}
      -\prescript{F}{}CLV_{i} \left(\delta, F_{i}, \boldsymbol{r}^{\text {acceptant }}_i, d_{i}, T\right)-c & \text { if } y_{i}=0\\
      \text{ }\\
      \gamma \cdot\Big[\prescript{F}{}CLV_{i}\left(p_{i}-\delta, F_{i}, \boldsymbol{r}^{\text {acceptant }}_i, d_{i}, T\right)& \text { if } y_{i}=1\\\phantom{ii} -\prescript{F}{}CLV_{i}\left(p_{i}, F_{i}, \boldsymbol{r}^{\text {lapser }}_i, d_{i}, T\right)\Big]-c .
    \end{cases}       
\end{equation}

That last equation can seem obscure at first glance. It simply assigns to each individual the expected profit or loss that would result from targeting her with a given lapse management strategy. A positive amount for subject $i$ means that targeting her would generate profit, whereas a negative one would lead to a loss for the insurer. We can take the example of a hypothetical scenario where $p_i=3\%$, $\delta = 0.05\%$, $\gamma=10\%$ and $c=10$ euros. It would generate $z_i$s taking values from $-234,614$€ to $53,066$€ with a mean of $-218$€ and a median of $-55$€. Different scenarios would result in very different distributions for the $z_i$'s.\\

Eventually, we define $\tilde{y}_{i}$ as a binary target variable indicating for policyholder $i$ if the individual expected retention gain resulting from a given retention strategy is a profit or a loss. More specifically, we define $\tilde{y}_{i}$ as
\begin{equation}
    \tilde{y}_{i}=\left\{\begin{array}{l}
    1 \text { if } z_i >0 \\
    0 \text { if } z_i \leq 0
    \end{array}.\right.
\end{equation}
\begin{mdframed}[hidealllines=true]
\begin{remark}
A subject in the dataset for which $y_i=0$ would produce $\tilde{y}_{i}=0$, whereas one for which $y_i=1$ could produce $\tilde{y}_{i}=0$ or $\tilde{y}_{i}=1$. In other words, it is never profitable for the insurer to offer an incentive to a subject that would not have lapsed. On the other hand, offering that same incentive to a lapser can be profitable. However, depending on the subject's features and the lapse management strategy parameters, it can also lead to a loss.
\end{remark}
\end{mdframed}

We can now include $\tilde{y}_{i}$ as a new binary target variable in our models and directly consider $RG$ as the global evaluation metric in the tree-based models we consider.
\\
We can now compare two models: the classical one with $y_i$ as a target variable and accuracy as the evaluation metric; and the CLV-augmented one with $\tilde{y}_{i}$ as a target variable and $RG$ as the evaluation metric.\\
\\
Intuitively, the former tries to predict whether a policyholder will lapse and tune its parameters by minimizing the misclassification rate. On the other hand, the latter aims at predicting whether applying a given retention strategy to the $i$-th individual will be profitable for the insurer and tune its parameters by maximizing the global expected retention gain.

\section{Methodology}\label{sec:Methodology}
In Section~\ref{sec:Framework}, we described a business-oriented framework, augmenting lapse management strategy with an evaluation metric based on the future CLV of every subject. Evaluating this metric requires computing $r^{\text{acceptant}}$ and $r^{\text{lapser}}$, the matrices of size $(n,T+1)$ containing for every subject, survival probabilities that we detail below. This individual survival analysis differs from \cite{loisel_piette_tsai_2021}'s work where $r^{\text{lapser}}$ is estimated globally and takes the same value for every policyholder regardless of their characteristics and where $r^{\text{acceptant}} = 1$ for any subject and at any time, ignoring the fact that an ``acceptant"'s policy can end with the policyholder's death.\\
Given this framework, we propose a two-step methodology: firstly, we detail how this survival analysis is carried out to model those retention parameters, and secondly, we explain how we use them for training tree-based classification models.
\subsection{Step 1: Modeling \texorpdfstring{$r^{\text{acceptant}}$}{r_acceptant} and \texorpdfstring{$r^{\text{lapser}}$}{r_lapser}}
We recall that a given subject's policy can end with lapse or death, and the policy is considered active if competing events are yet to occur. Furthermore, while a lapser's policy can end with lapse or death, whatever comes first, an acceptant one can only end with death.\\
$r^{\text{lapser}}$ represents the probability that the policy of subject $i$ is still active at time $t$, given that the subject is labeled as a lapser - $\text{EVENT}=1$ - at $t=0$. Predicting these overall conditional survival probabilities with competing risks can be achieved by creating a combined outcome: the policy ends with death or lapse, whichever comes first. To compute $r^{\text{lapser}}$ in practice, we recode the competing events as a combined event. This approach is compatible with any survival analysis method regarding statistical guarantees.\\
\\
Conversely, $r^{\text{acceptant}}$ represents the probability that the policy of subject $i$ is still active at time $t$, given that the subject is not labeled as a lapser - $\text{EVENT}=0$ or $2$ - at $t=0$. This estimation is more complex as we must dissociate the risks of lapse and death. These causes being mutually exclusive, a competing risks methodology is well-suited to estimate $r^{\text{acceptant}}$ \cite{Laurent2016-qt}.\\
\\
It is also important to note that here, $r^{\text{lapser}}$ is modeled on subjects that have lapsed in the past - they may have been offered an incentive in the past, this is unknown - and not on subjects that have been offered an incentive that they declined. Our framework makes the implicit hypothesis that both behaviors are alike. It is more intuitive for $r^{\text{acceptant}}$ as a subject that has not lapsed in the past would have accepted any incentive if offered.

\subsubsection{Competing risks frameworks}
\label{com_risk_fram}
We are aware that improvements of our model over \cite{loisel_piette_tsai_2021}'s approach, require the analysis of both the risks of lapse and death, thus a competing risk setting. As detailed in Appendix~\ref{appendix_competing_risks}, several regression models exist to estimate the global hazard and the hazard of one risk in such settings: cause-specific and subdistribution models. They account for competing risks differently, obtaining different hazard functions and thus have distinct advantages, drawbacks and interpretations. These differences are discussed in \cite{milhaud2018}, where the authors also considered a competing risk framework for lapse prediction.\\
After discussions detailed in Appendix~\ref{appendix_competing_risks}, the simplicity of a cause-specific approach and the fact that it can be adapted to any survival method, including tree-based ones, oriented our choice towards it. We then computed  $r^{\text{acceptant}}$ and $r^{\text{lapser}}$ with three different methods - Cox model, Random Survival Forest and Gradient Boosting Survival Model - and retained the best one. These methods are shortly described in the following sections.

\subsubsection{Cox proportional hazard model}
One of the most common survival models is the Cox proportional hazard (CPH) model (\cite{cox1972a}). It postulates that the hazard function can be modelized as the product of a time-dependent and a covariate-dependent functions.
The hazard function at time $t$ for subject $i$ with covariate vector $\boldsymbol{X_i}$, under Cox proportionnal hazard model can be expressed as\\
$$
\underbrace{\lambda(t | X_{i}^{1}, X_{i}^{2}, \dots)}_{\text{hazard function}} = \lambda(t | \boldsymbol{X_i}) = \overbrace{\lambda_0(t)}^{\text{baseline hazard}} \underbrace{e^{\overbrace{\left( \boldsymbol{X_i} \cdot \beta_i \right)}^{\text{log-partial hazard}}}}_ {\text{partial hazard}} .
$$
It is crucial to note that in this model, the hazard function is the product of the baseline hazard, which only varies with time, and the partial hazard, which only varies depending on the covariates.
The parameters of this model are the $\beta$, and they can easily be estimated with a maximum likelihood approach. Their estimation can be carried out without having to model $\lambda_{0}(t)$ - which is why CPH is considered semi-parametric.\\
We use Python and lifelines (\cite{Davidson-Pilon2019}) to implement it. We specify a spline estimation for the baseline hazard function. We select the covariates and model parameters using AIC (\cite{akaike-a}) and use the concordance index (\cite{harrell1982evaluating}) to compare CPH to other models. The concordance index - or Harrel's c-index or simply c-index - is a metric to evaluate the predictions made by a survival model. It can be interpreted as a generalization of the area under a receiver operating characteristic (ROC) curve (\cite{hanley1982a}) - or AUC - in a survival setting with censored data.

\subsubsection{Random Survival Forest}
Survival trees have been extensively studied for a long time, and a complete review of such existing methods up to 2011 can be found in \cite{BouHamadreview}. The most important thing to understand is that a survival tree can be created by modifying the splitting criterion of a regular tree. Most survival tree algorithms are designed with a split function that aims to maximize the separation of the resulting child nodes in terms of survival profiles. This separation between nodes is estimated by maximizing the log-rank statistic (\cite{mantel1966a, LeBlanc.1993}). Each terminal node of a survival tree contains a survival profile from which we can derive the survival and cumulative hazard function.\\
\\
An RSF is the counterpart of a random forest (see Appendix~\ref{appendix_rsf}) for such survival trees. It has been developed in \cite{RSFIshwaran08} and extended for competing risks a few years after (\cite{ishwaran2014a}). A prediction with RSF for a given subject is made by getting his/her survival profile in each tree in the forest. His/her corresponding survival and cumulative hazard function are estimated in each tree with Kaplan-Meier and Nelson-Aalen estimators, respectively. Eventually, the aggregation of those single-tree estimates constitutes the RSF's prediction. \\
\\
We use Python and sksurv (\cite{sksurv}) to implement RSF, and we tune and evaluate our model using the concordance index.
\begin{mdframed}[hidealllines=true]
\begin{remark}\label{remarkRSF}
Sksurv allows us to use RSF with a cause-specific consideration of the competing risks. To this day, sksurv does not have a subdistribution competing risks model, whereas its R implementation \emph{randomForestSRC} does (\cite{RSFR}).\\
Moreover, a severe limitation of that approach is that predictions can only be made at time points observed in the training set. Concretely, this prevents us from using RSF to extrapolate survival and hazard functions to unobserved time points.
\end{remark}
\end{mdframed}

\subsubsection{Gradient Boosting Survival Model}

In the same way Random Forest has a survival counterpart, this is also true for Gradient Boosting approaches. An essential distinction between classical boosting algorithms (see Appendix~\ref{appendix_gbsm}) and Gradient Boosting Survival Model (GBSM) lies in its loss function. The loss function that we use with GBSM is the partial likelihood loss of a CPH model, and the optimization in such a model is achieved by maximizing a slightly modified log-partial likelihood function, 
\begin{equation*}
  \arg \min_{f} \quad \sum_{i=1}^n \delta_i \left[ f(\boldsymbol{X_i})
- \log \left( \sum_{j \in {g}_i} e^{(f(\boldsymbol{X_j}))} \right) \right], 
\end{equation*}
where $\delta_i$ is the event indicator and $f(\boldsymbol{X_i})$ is GBSM's prediction for the $i$-th subject, with a covariate vector $\boldsymbol{X_i}$. ${g}_i$ is the tree leaf including subject $i$.
\\
Similarly to RSF, we use Python and sksurv (\cite{sksurv}) to implement GBSM. We tune and evaluate our model using the concordance index. Remark~\ref{remarkRSF} also applies here.

\subsubsection{Final modeling choice}
Our analysis shows that, based on concordance index, RSF and GSBM both outperformed a semiparametric Cox model in our study case. Regarding interpretability, we note that the feature importance analysis is very similar between the three models. All the details about the final concordance index scores, covariates importance and various plot for further analysis are available in Appendix~\ref{appendix_surv}.\\
\\
In the following sections, we decide to retain \textbf{GBSM} for the modeling of $r^{\text{acceptant}}$ and $r^{\text{lapser}}$ as it has the best concordance index.
\vspace{3 mm}
\begin{mdframed}[hidealllines=true]
\begin{remark}
As this study aims to be business-oriented and favor real-life decision-making, it is crucial to note that the computation times for fitting these different models are very different and could potentially be a huge constraint for real operational deployment. Specific computation times differ greatly depending on various factors, such as the number of subjects or features considered, the computation power or parallelization ability at disposal, for instance. However, we can still give here an order of magnitude for those differences. If the tuning and fitting process for CPH can last a few tens of seconds, it lasts hours for RSF and tens of hours for GBSM.
\end{remark}
\end{mdframed}



\subsection{Step 2: Classification tasks}\label{classif_step}
Our work focuses on lapse management with tree-based models. It aims to answer the question: which policyholders would be worth targeting with a lapse management strategy to maximize the expected T-year profit for the insurer? We will consider a single tree built with Breiman's CART algorithm, Random Forest, XGBoost, and RE-EM trees. The following sections detail how those different approaches work. Those models will be compared on two different classification tasks; and tuned with two different evaluation metrics, a statistical metric and a business-related one.

\paragraph{\underline{On $y_i$}}\label{y}
First, we will use a classical lapse prediction framework to model the policyholder's behavior. Each policyholder will be labeled as a lapser or a non-lapser with a binary outcome $y_{i}$. Our first batch of models will be trained with $y_{i}$ as a response variable and produce predictions $\hat{y}_{i}$. $\operatorname{Accuracy}(y, \hat{y})$, which is undoubtedly the most intuitive performance measure for binary classification, is defined as the proportion of correctly predicted observations over all observations. It is widely used for churn analysis and appears to be a satisfying performance measure for relatively balanced outcomes - 22\% of all observed subjects being lapsers - in binary classification problems. We will use it as an evaluation metric in a 10-fold cross-validation step for tuning our models.\\
We know that more advanced evaluation metrics are available for binary classification, including the recall, the $\operatorname{F_{\beta}}$ score family (\cite{Chinchor1992}), the AUC under the ROC or under the Precision-Recall curve, the Brier Score (\cite{BRIER1950}) and lift curve. They are standard evaluation metrics in classification and provide valuable insights into the model's performance, they are also frequently used in the applied binary classification literature, especially in the presence of a significant imbalance in the data (\cite{Haibo2009}). However, in this paper, the mildness of the imbalance of $y_i$ and our will to compare a customer-centered framework to representative real-world practices encourages us to use accuracy as a comparison. One of the goals of this article is to demonstrate that some of the current practices in real-world applications, based on statistical metrics such as accuracy can be significantly improved by considering a profit-driven target variable and evaluation metric. We are aware that accuracy may not be an optimal choice of evaluation metric for binary prediction in general and churn or lapse analysis specifically, but it seems representative of what practitioners use (see Table 2 from \cite{duchemin2021} for example), as it is suggested in \cite{loisel_piette_tsai_2021}. We do not aim at comparing our framework against the best existing methods but rather against the most representative. Nevertheless, the numerical results of Table~\ref{tab:Numerical_results_1} have also been obtained with recall, F1-score, and AUC for tuning and cross-validation and some are available in Appendix~\ref{other_metrics}: the conclusions obtained with such measures are similar to those obtained with accuracy. Thus, in this article and as in \cite{loisel_piette_tsai_2021}, we will only select, evaluate and discuss the models in the light of accuracy.\\

\paragraph{\underline{On $\tilde{y}_i$}}\label{y_tilde}
Secondly, we will use the CLV-Augmented lapse prediction framework, detailed in Section~\ref{sec:Framework}. Each policyholder will be labeled as a targeted lapser or a non-targeted policyholder with the binary outcome $\tilde{y}_i$ and prediction for that outcome are denoted $\hat{\tilde{y_i}}$. \\
\begin{mdframed}[hidealllines=true]
\begin{remark}
Note that whenever $y_{i}=0$, we also have $\tilde{y}_i=0$. In other words, if subject $i$ does not intend to lapse, it is never worth proposing her an incentive: the subject will accept it with probability $1$ and would not have lapsed.\\
On the other hand, when $y_{i}=1$, it corresponds to either $\tilde{y}_i=1$ or $\tilde{y}_i=0$. In other words, if subject $i$ is labeled as a lapser, it does not necessarily mean it is worth targeting her. From the insurer's point of view, some policies are better off lapsed. $\tilde{y}$ can be seen as a more detailed version of $y_{i}$ as it carries not only behavioral information regarding lapse but also a profitability one.
\end{remark}
\end{mdframed}

We thus train a second batch of models with $\tilde{y}_i$ as a response variable. We use $RG$ as an evaluation metric in a 10-fold cross-validation step for tuning these models.
\begin{mdframed}[
  leftmargin=0.5\parindent,
  rightmargin=0.5\parindent,
  skipabove=\topsep,
  skipbelow=\topsep,
  hidealllines=true
  ]
  \emph{Summary of our methodology:} First, we train a CART, RF and XGBoost models with $y_{i}$ as a binary target variable and accuracy as a tuning evaluation metric. \\
  Then we train them with $\tilde{y}_i$ as a binary target variable and $RG$ as a tuning evaluation metric.\\
  Finally, we train and test all six models on different random samples of our dataset and keep track of the model's classification performance on all of them and for various retention strategies for comparison's sake.
\end{mdframed}
\vspace{2 mm}
The sections below briefly introduce the tree-based model we selected before displaying how they performed in various lapse management scenarios.

\subsubsection{CART}
CART (\textit{Classification And Regression Trees}) is an algorithm developed by \cite{breiman1984classification} that consists of recursively partitioning the covariate space. It is a widespread, intuitive and flexible algorithm that handles regression and classification problems.

\subsubsection{Random forest}
A natural idea to correct CART's instability and enhance its prediction accuracy is the aggregation of a significant number of single trees, each grown on different subsamples of the dataset. A random forest (RF by \cite{breiman2001random}) is a tree-based bagging procedure where each tree is grown on randomly drawn observations and contains splits considering only randomly drawn covariates. 

\subsubsection{XGBoost}
Other tree-based approaches have been designed to reduce the instability of a single-tree model. Model boosting is an adaptative technique, first developed by Freund et Shapire (\cite{Freund1996}), that does not rely on the aggregation of independent weaker models but rather on the aggregation of weak models built sequentially, one after the other. XGBoost (\cite{Chen2016}) is a widespread and performant tree-boosting model that relies on a gradient-boosting step and provides a very optimized parallelized procedure. It is considered a state-of-the-art library for various prediction problems.\\
\\
The interested reader can find more detailed explanations about CART, RF and XGBoost mechanisms in the aforementioned references. For these modeling approaches, we used Python and sklearn (\cite{pedregosa2011scikit}).


\section{Real-life application}\label{results}
Based on the real life-insurance dataset at our disposal (described in Section~\ref{sec:Data}), we use the survival model we selected and estimate $r^{\text{acceptant}}$ and $r^{\text{lapser}}$ for every individual. This allows us to compute the individual CLVs, RGs, $z_i$'s and $\tilde{y}_{i}$. We have already defined what a strategy is (see Definition~\ref{LMS}), and we can thus apply our classification methodology to various retention strategies.

\subsection{Considered lapse management strategies}
The strategies considered are based on several criteria. First, we selected realistic strategy parameters and time horizons based on actual retention campaigns led by life insurers. Moreover, we chose to present strategies that illustrate the exhaustive list of conclusions and discussions that are carried out in the next section. Finally, we also incorporated strategies that are ``obviously bad" in the sense that such strategies would necessarily lead to a loss for the insurer. Such extreme scenarios will supplement our discussions. In any case, we consider $p_{i}$ and $d_{i}$ to be constant in our application, as both those parameters were not estimated at the individual level by the life insurer that provided with the dataset.\\
Results related to the 64 considered LMS are given in Appendix~\ref{more_lms_results}. Our analysis showed that all considered LMS results can be split into 5 categories depending on how applying our framework impacted their expected retention gain over a naïve targeting. We have realistic profitable strategies that are improved by our framework, but also highly loss-inducing, moderate loss-inducing, highly profitable and unrealistically highly profitable strategies. We refer to the LMS displayed in Table~\ref{tab:LMS} as representative strategies as they all belong to one of those categories. Numerical results regarding the most representative strategies can be found in Section~\ref{Results_section} and related comments on how to read these tables are given in Section~\ref{Comment_section}.\\
\\
\begin{adjustbox}{width=.5\columnwidth,center}
    \begin{tabular}{c|cccccc}
    \multicolumn{1}{c}{Scenarios} & p     & $\delta$ & $\gamma$ & c     & d     & T \\
    \midrule
    A-1     & 2.50\% & 0.04\% & 25\%  & 10    & 1.50\% &                     5  \\
    A-5     & 2.50\% & 0.04\% & 5\%   & 10    & 1.50\% &                     5  \\
    A-25    & 5.00\% & 0.10\% & 25\%  & 10    & 1.50\% &                     5  \\
    B-6    & 2.50\% & 0.08\% & 10\%  & 10    & 1.50\% &                   20  \\
    B-27    & 5.00\% & 0.20\% & 20\%  & 100   & 1.50\% &                     5  \\
    \end{tabular}%
\end{adjustbox}
\captionof{table}{Insightful LMS} 
\label{tab:LMS}
\subsection{Numerical results}\label{Results_section}
\begin{adjustbox}{width=1\textwidth}
    \begin{tabular}{cclcrrrrrrr}
    \multirow{2}[1]{*}{N°} & \multirow{2}[1]{*}{time (s)} & \multicolumn{1}{c}{\multirow{2}[1]{*}{Model}} & \multicolumn{1}{c}{\multirow{2}[1]{*}{\% target diff (\% of 1's)}} & \multicolumn{2}{c}{Accuracy} & \multicolumn{2}{c}{ Retention gain } & \multicolumn{2}{c}{ RG/target } & \multicolumn{1}{c}{\multirow{2}[1]{*}{Improvement}} \\
          &       &       &       & \multicolumn{1}{c}{$y_i$} & \multicolumn{1}{c}{$\tilde{y}_i$} & \multicolumn{1}{c}{ $y_i$ } & \multicolumn{1}{c}{ $\tilde{y}_i$ } & \multicolumn{1}{c}{ $y_i$ } & \multicolumn{1}{c}{ $\tilde{y}_i$ } &  \\
    \midrule
    \multirow{3}[2]{*}{A-1} & \multirow{3}[2]{*}{4949} & CART  & \multirow{3}[2]{*}{62.6\% (8.2\%)} & 92.3\% & 85.3\% &                114 661  &              219 655  &               4.48  &            38.20  & 91.6\% \\
          &       & RF    &       & 92.9\% & 85.4\% &              232 314  &             287 884  &               9.82  &            56.65  & 23.90\% \\
          &       & XGB   &       & 93.4\% & 85.8\% &             243 365  &             324 952  &                9.61  &            54.64  & 33.50\% \\
    \midrule
    \multirow{3}[2]{*}{A-5} & \multirow{3}[2]{*}{4753} & CART  & \multirow{3}[2]{*}{86.7\% (2.9\%)} & 92.3\% & 83.6\% & -          514 477  & -           112 372  & -        20.08  & -        86.48  & 78.20\% \\
          &       & RF    &       & 92.9\% & 83.4\% & -         323 544  & -               3 937  & -         13.65  & -        28.28  & 98.80\% \\
          &       & XGB   &       & 93.4\% & 83.3\% & -         383 004  &                           0    & -          15.14  &                     0    & 100.00\% \\
    \midrule
    \multirow{3}[2]{*}{A-25} & \multirow{3}[2]{*}{5379} & CART  & \multirow{3}[2]{*}{31.0\% (15.2\%)} & 92.3\% & 89.2\% &          4 160 423  &         3 882 623  &          162.44  &          241.06  & -6.70\% \\
          &       & RF    &       & 92.9\% & 89.5\% &          4 018 432  &          3 666 219  &          169.65  &         249.54  & -8.80\% \\
          &       & XGB   &       & 93.4\% & 90.0\% &          4 455 108  &          4 410 629  &          176.09  &         267.87  & -1.00\% \\
    \midrule
    \multirow{3}[2]{*}{B-6} & \multirow{3}[2]{*}{5906} & CART  & \multirow{3}[2]{*}{36.6\% (13.9\%)} & 92.3\% & 88.8\% &              705 721  &             922 490  &            27.69  &             60.21  & 30.70\% \\
          &       & RF    &       & 92.9\% & 88.9\% &           1 352 182  &          1 269 349  &              57.11  &            97.63  & -6.10\% \\
          &       & XGB   &       & 93.4\% & 89.6\% &          1 342 882  &          1 428 722  &            53.09  &            96.76  & 6.40\% \\
    \midrule
     \multirow{3}[2]{*}{B-27} & \multirow{3}[2]{*}{4811} & CART  & \multirow{3}[2]{*}{73.9\% (5.7\%)} & 92.3\% & 84.3\% & -         694 436  &              751 404  & -         27.16  &         226.99  & 208.20\% \\
          &       & RF    &       & 92.9\% & 84.4\% & -              13 512  &           1 018 369  & -            0.51  &         356.48  & 7637.00\% \\
          &       & XGB   &       & 93.4\% & 84.7\% & -            38 050  &          1 253 252  & -            1.55  &         345.94  & 3393.70\% \\
    \bottomrule
    \end{tabular}%
\end{adjustbox}
\captionof{table}{Means of the results obtained on considered LMS} 
\label{tab:Numerical_results_1}

\subsection{Comments}\label{Comment_section}
Several terms in the two previous tables need to be explained. ``\% target diff" represents how different $y$ and $\tilde{y}$ are. It is the percentage of subjects for which  $y_i = 1$ and $\tilde{y}_i = 0$: in other words, the proportion of lapsers not worth targeting with a given strategy. The quantity ``\% of 1's" represents the proportion of ones in $\tilde{y}$ the target variable. It is to be compared with the 22\% of ones in $y$: the proposed framework's imbalance increases with ``\% target diff".\\
Then the table shows the 10-fold cross-validated mean accuracies, retention gains and RG/target with two methodologies: the columns denoted $y_i$ represent the metrics obtained by a model with $y_i$ as a response variable and accuracy as an evaluation metric, and the columns denoted $\tilde{y}_i$ represent the metric obtained by a model with $\tilde{y}_i$ as a response variable and $RG$ as an evaluation metric.\\
RG/target represents the achieved retention gain for every targeted individual, for $y_i$, it is $RG/\sum_i \hat{y_i}$, for $\tilde{y}_i$ it is $RG/\sum_i \hat{\tilde{y_i}}$. Eventually, ``Improvement" represents the percentage of improvement between the $RG$ obtained with a classification on $y_i$ and the gain obtained with a classification on $\tilde{y}_i$. As the reported financial information was distorted for confidentiality reasons (see Section~\ref{sec:Data}), relative measures such as ``Improvement" are certainly more informative than absolute ones such as $RG$.\\
Some LMS are worth focusing on. For every strategy, we display its 10-fold cross-validated results: 10\% of the dataset acting as an out-of-sample validation set at every fold. Every model is tuned by cross-validation within every fold. The boxplots below summarize some typical key results illustrated by several strategies. Those results will be discussed in Section~\ref{Discussion}.\\

\begin{minipage}{1\linewidth}
\makebox[\linewidth]{

  \includegraphics[width=0.65\linewidth]{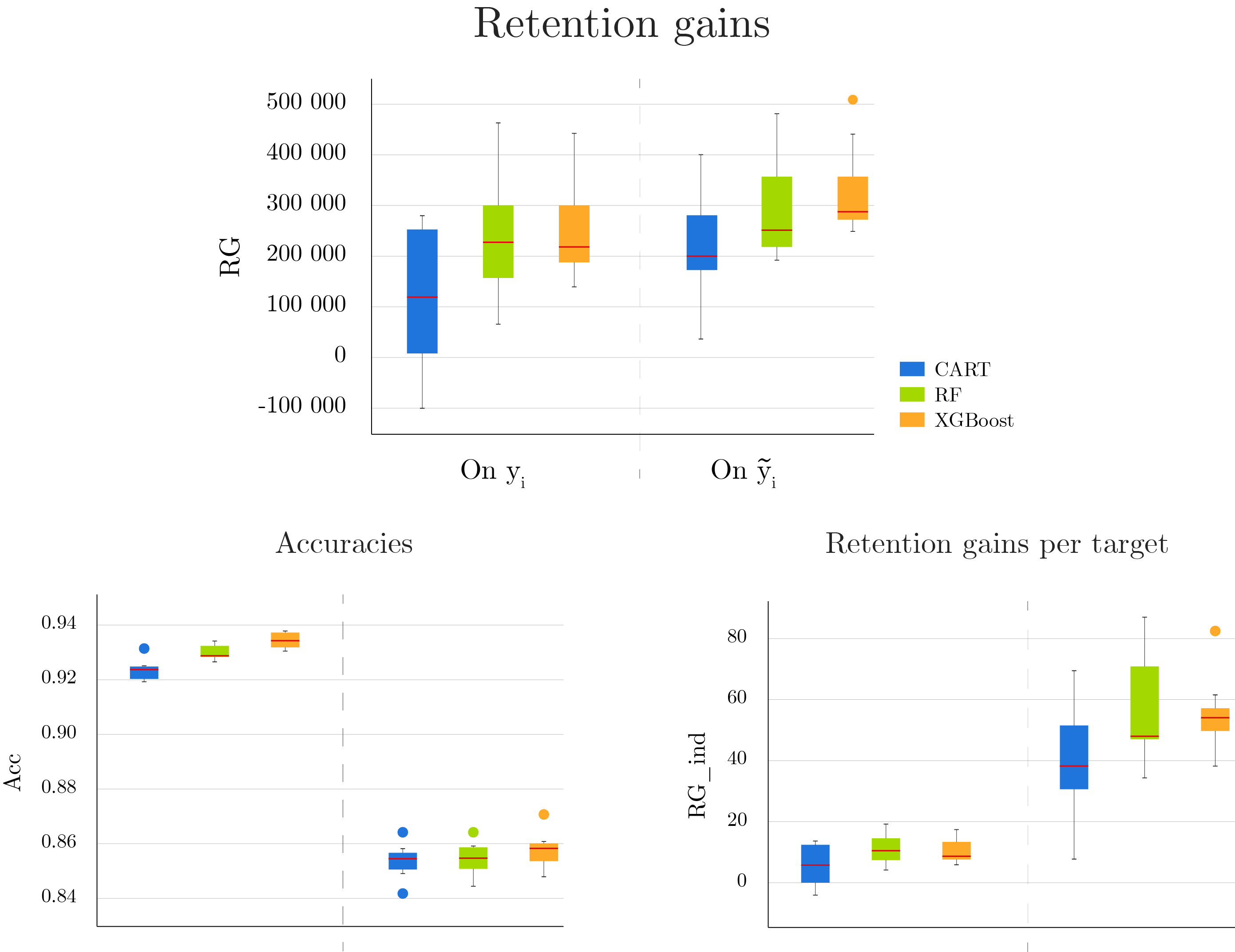}}
\captionof{figure}{Strategy n°A-1: (Positive result on $y_i$ and an improved result on $\tilde{y}_i$.)}\label{lms1_graphs}
\end{minipage}
\begin{minipage}{1\linewidth}
\makebox[\linewidth]{
  \includegraphics[width=0.65\linewidth]{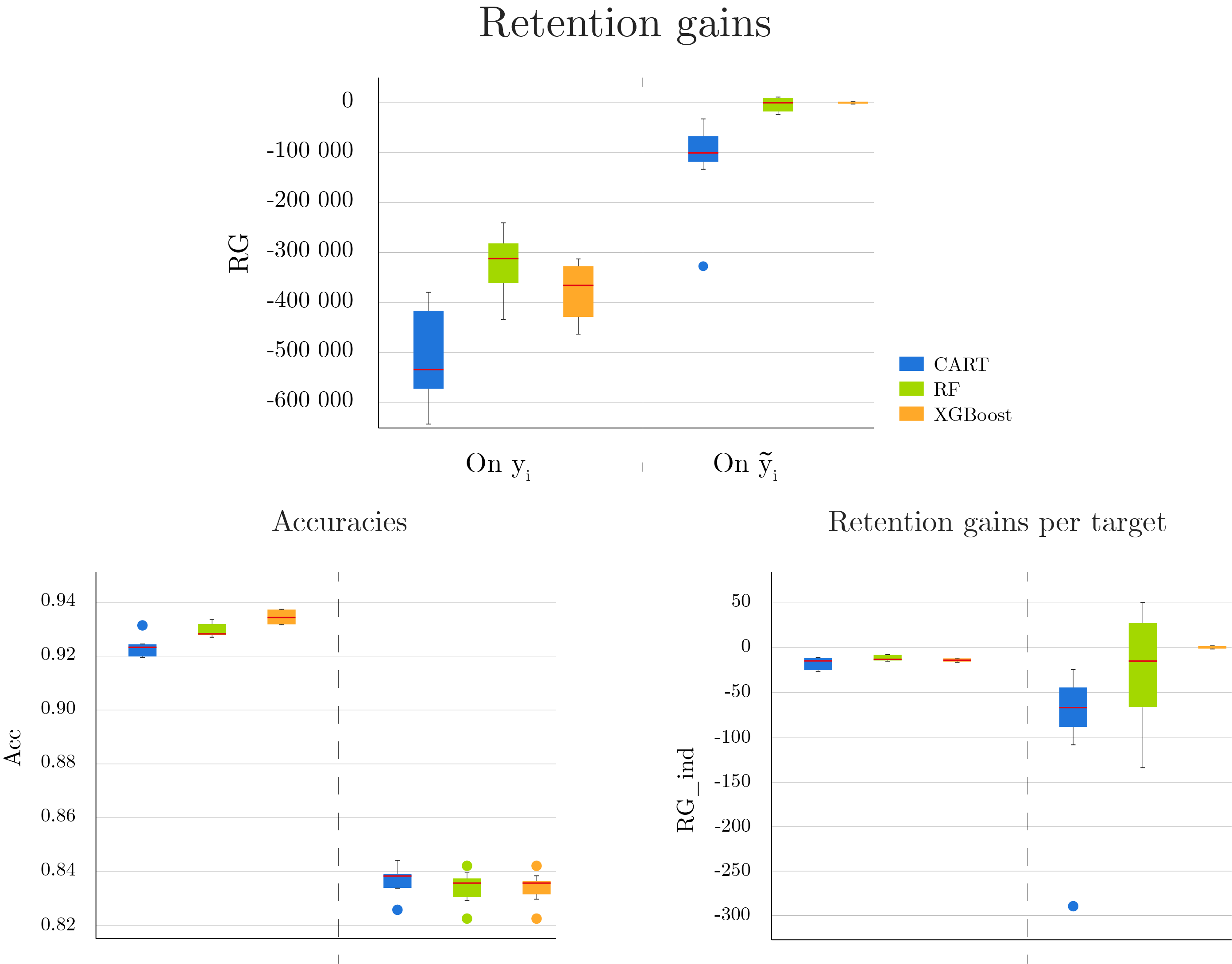}}
\captionof{figure}{Strategy n°A-5: (Very negative result on $y_i$ and a loss-limiting result on $\tilde{y}_i$.)}\label{lms4_graphs}
\end{minipage}
\begin{minipage}{1\linewidth}
\makebox[\linewidth]{
  \includegraphics[width=0.65\linewidth]{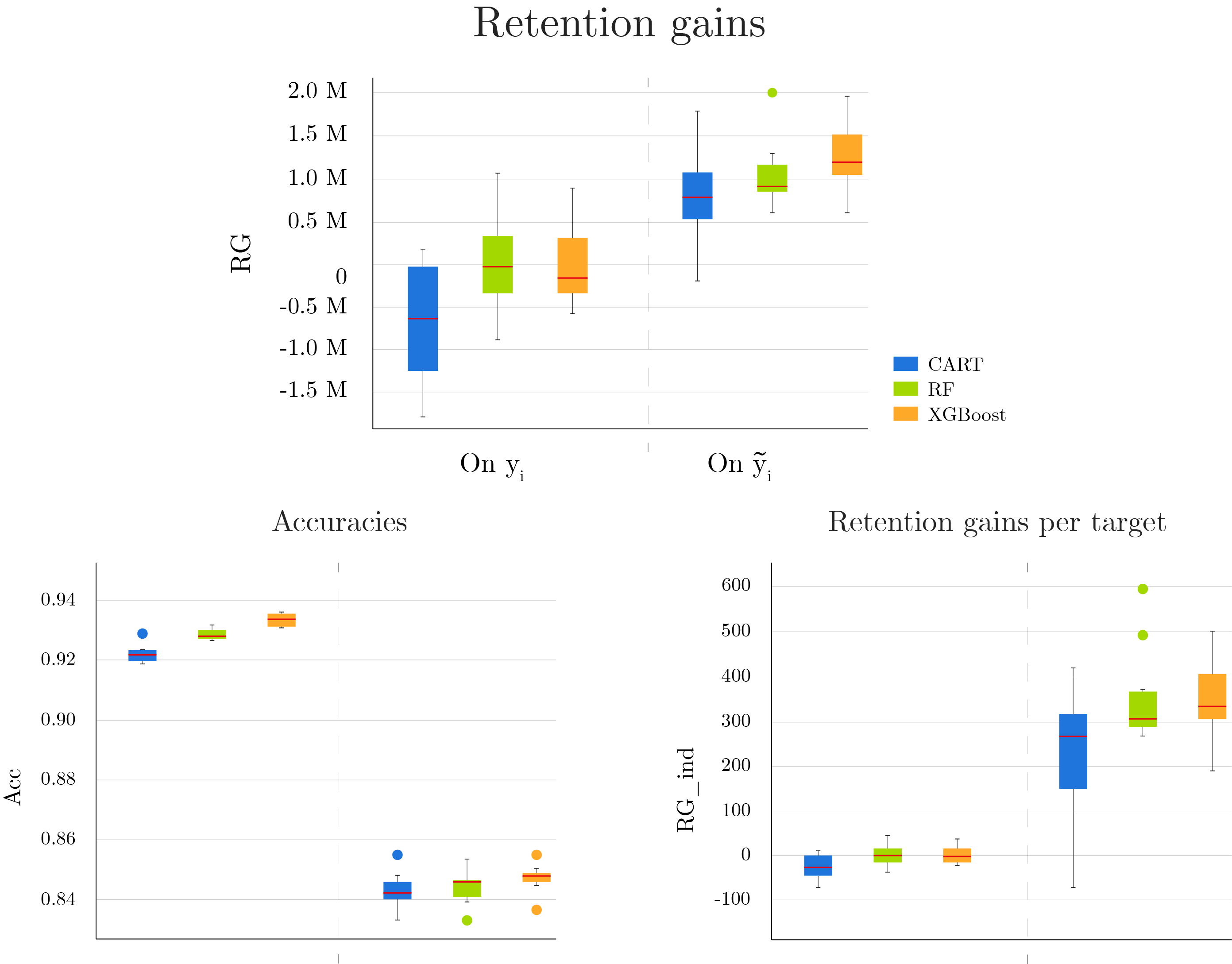}}
\captionof{figure}{Strategy n°B-27: (Negative result on $y_i$ and positive one on $\tilde{y}_i$)}\label{lms30_graphs}
\end{minipage}
\begin{minipage}{1\linewidth}
\makebox[\linewidth]{
  \includegraphics[width=0.65\linewidth]{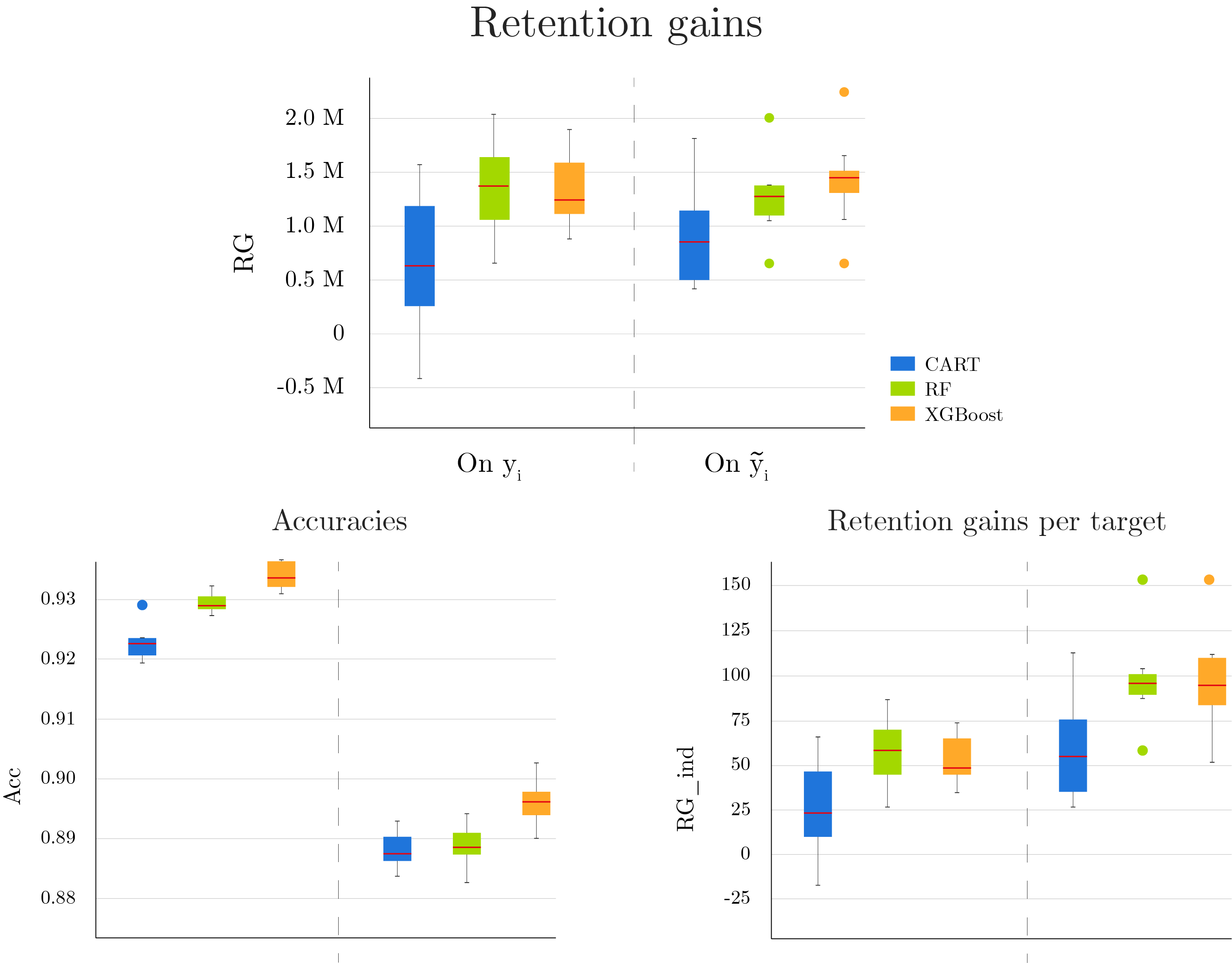}}
\captionof{figure}{Strategy n°B-6: (High positive result on $y_i$ slightly improved with $\tilde{y}_i$.)}\label{lms21_graphs}
\end{minipage}
\begin{minipage}{1\linewidth}
\makebox[\linewidth]{
  \includegraphics[width=0.65\linewidth]{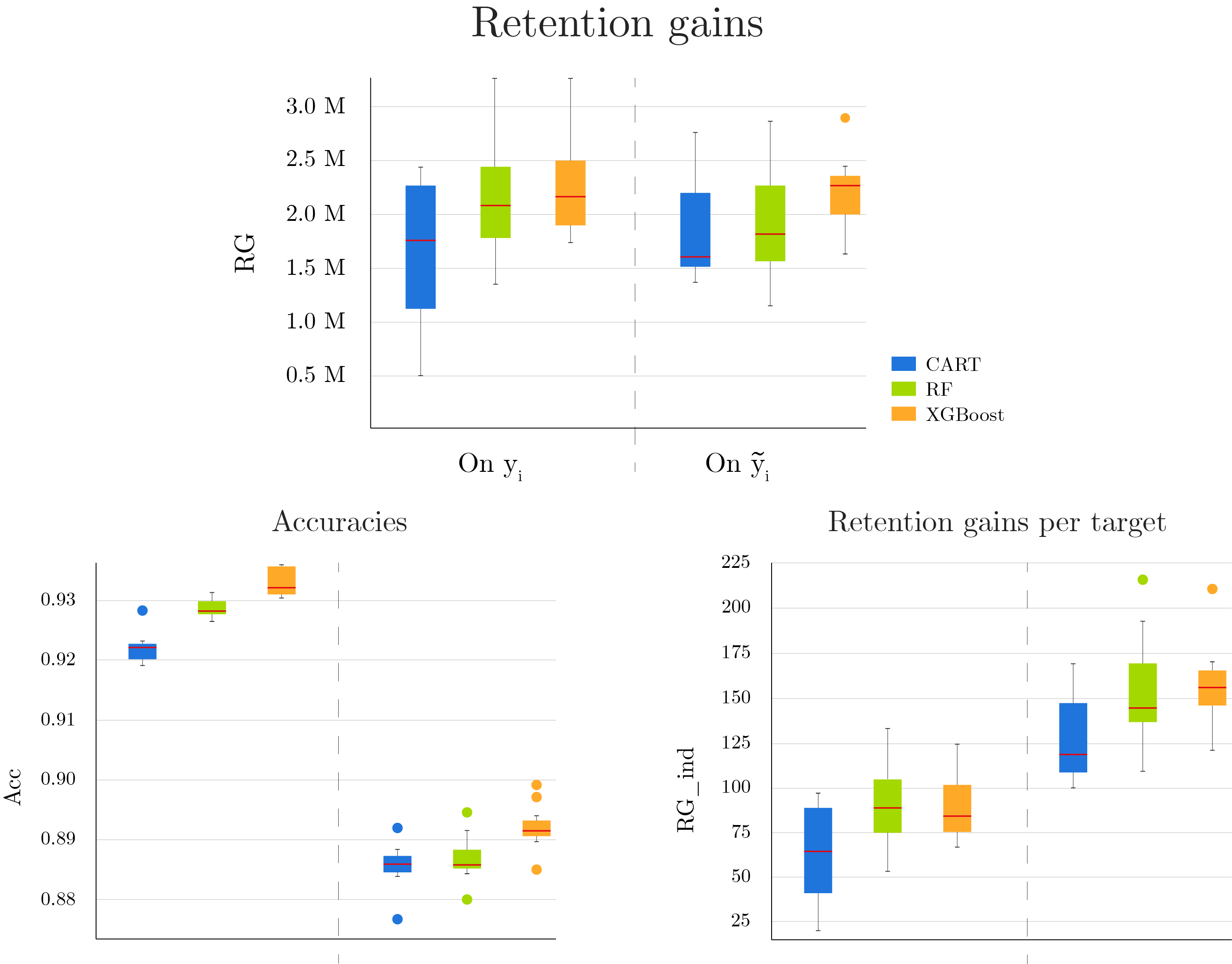}}
\captionof{figure}{Strategy n°A-25: (Results on $y_i$ better than results on $\tilde{y}_i$.)}\label{lms13_graphs}
\end{minipage}
\begin{mdframed}[hidealllines=true]
\begin{remark}
With considerable computation power and great parallelization, the results for all strategies - see other strategies in Appendix~\ref{more_lms_results} - were obtained with a wall time of less than 4 days and a CPU time of more than 100 days.
\end{remark}
\end{mdframed}

\section{Discussion}\label{Discussion}
\subsection{General statements}\label{General}
As expected and shown in the actuarial literature, RF and XGBoost perform globally better than CART regarding mean accuracy and RG. It is true for all LMS considered in Table~\ref{more_lms_considered2}. Globally, XGBoost is more consistent and is the best model in most scenarios, both with and without the CLV-based measure. It is only outperformed by RF in strategies n°A-7, A-11, A-14, A-29, B-7, B-14 and B-31.\\
As expected, by design, the vast majority of strategies, including all the realistic ones, show that a classification on $\tilde{y}_i$ produces a targeting that yields better RG than a classification on ${y}_i$. Conversely, a classification on $y_i$ produces a targeting that delivers better accuracies regarding whether a policyholder will churn than a classification on $\tilde{y}_i$. These results were expected because of the models' respective objectives. Even if it is not surprising, it once again shows that for an insurer, lapse prediction and lapse management strategy are two very different prediction problems, often treated as similar ones.\\
\\
Our CLV-augmented model shows different behavior depending on the strategy considered. As highlighted by Figure~\ref{lms1_graphs}, a model on $y_i$ is greatly improved by our framework regarding RG and RG/target. Conversely, its accuracy in lapse prediction is not optimal. \\
An attractive property of our framework can be observed in Figure~\ref{lms4_graphs}: it yields loss-limiting targeting. When the LMS considered is too aggressive, it will usually prefer to predict that an LMS should not be applied at all ($\forall i, \hat{\tilde{y_i}}=0$), thus generating a $RG$ around 0€. This is made evident in some extreme strategies (LMS n°A-5, A-15, B-11, B-13 and B-16) and explains the presence of 0's in Table~\ref{tab:Numerical_results_1}.\\
On less extreme strategies, it shows to yield substantial improvement when classification on ${y}_i$ gives negative RG. That observation confirms what was already pointed out by \cite{loisel_piette_tsai_2021}: it can even turn a negative $RG$ into a positive one (see LMS n°A-8, B-8, B-12, B-23 and B-27 (Figure~\ref{lms30_graphs})) .\\
Our framework also improves a strategy where a classification on ${y}_i$ gives high RG. However, the improvement decreases as the difference between the total number of lapsers and the number of lapsers that would be profitable if retained is sizeable. An example of that is shown in Figure~\ref{lms21_graphs}.\\
\\
Finally, we can generate LMS for which our framework does not improve the expected RG. It is the case in LMS n°A-13, A-18 or A-27 (See e.g  Figure~\ref{lms13_graphs}). In LMS n°A-13, we can see that the mean of the RG is not improved, but the median is. In all those cases, the RG per target produced by the CLV-augmented model is greatly improved, indicating that a CLV-augmented strategy prefers to target fewer policyholders but only those who would generate high future profits. This last observation explains why a CLV-augmented LMS generates higher RGs when the cost of contact c is considerable. Indeed, the more costly a contact is, the more precise and specific a targeting strategy should be. \\
\\
Generally, we can collate the results of various LMS - excluding LMS n°B-27 that has a very high improvement ratio - to obtain a mean performance of our framework.

The average observed RG improvement of a CLV-augmented framework over the classical lapse one is 57,9\%\footnote{Using XGBoost}. If we weigh these results by the expected RGs, the average improvement is still 31,7\%. As a comparative result, it is reported in Section 6.2 of Loisel et al.'s work (\cite{loisel_piette_tsai_2021})  that they obtain improvements over that same classical framework between 18\% and 26\%, depending on the considered strategies. This emphasizes that by extending their work, we seem to improve on their results. Obviously, as we were not able to compare our results on the same data and strategies, and because our definitions of $RG$ differ, such a conclusion is to be treated cautiously.
\vspace{2 mm}
\subsection{Marketing decision making}
We already pointed out that the improvement of a lapse management strategy including CLV grows with the proportion of lapsers with a negative CLV (see Appendix~\ref{other_results}). Models resulting from our framework do not consider them as good targets. In fact, there is a Pearson correlation coefficient of 77\% between RG improvement and the proportion of target differences among the LMS detailed in Table~\ref{more_lms_considered2}. Of course, as the improvement ratio has no clear interpretation in some cases, this analysis should be carried out in more depth, separating the cases where both RG - with and without the inclusion of CLV - are positive from the cases where one of them is negative. By doing so, we observe that the Pearson correlation coefficient for LMS yielding positive RG regardless of the inclusion of CLV is even higher: 83\%.\\

In terms of targeting, it seems crucial to understand what differentiates a subject for which $y_i = 1$ and $\tilde{y}_i = 0$ from the others. An investigation of such policyholder profiles can be carried out for every lapse management strategy. We take the example of LMS n°A-1, where 62,6\% of policyholders were in that case (see Section~\ref{Results_section}). With that strategy, the profile of non-targeted lapsers indicates that
\begin{itemize}
    \item 57.2\% of them are men, similar to the entire dataset,
    \item 76.4\% of them contracted product n°1 whereas 72\% of all policyholders chose it,
    \item the mean seniority of their policy is 10.4 years compared to the 13.4 years for the complete dataset,
    \item the mean face amount of such policies is 12,156, whereas the average face amount for all considered policies is 40,263.
\end{itemize}

In that strategy, our framework indicates that marketing efforts on low seniority policyholders with low face amount policies are inefficient. Of course, this conclusion is only valid for the considered LMS; however, our framework allows us to conduct such analysis for any LMS and interpret the results at an individualized level. 

\subsection{Management rules decision making}
Sensitivity analysis of those results can highly benefit management rules decision-making. This framework serves as a tool that compares future hypothetical lapse management strategies in order to choose the best one - among realistic scenarios -. It can also be used to tune a given strategy by answering questions like:

\begin{itemize}
    \item For which incentive $\delta$ the retention strategy becomes profitable ?
    \item For which acceptance probability $\gamma$ the retention strategy becomes profitable ?
    \item With a given budget, what is the optimal list of policies that should be targeted? 
    \item At which horizon $T$, the retention strategy become profitable ? In other words, when can the insurer expect a return on investment?
\end{itemize}

Answering these questions constitutes a 1-parameter sensitivity analysis. In our framework, six parameters influence the expected retention gain $(p, \delta, \gamma, c, d, T)$. \\
We can argue that among them are three structural parameters that are insurer's dependent and not linked to the external state of the world: $\delta$, $\gamma$ and $c$. Among them, the contact cost $c$ is more or less fixed and can not be easily changed by the insurer. Conversely, $\delta$ and $\gamma$ are to be chosen by the insurer. Moreover, they also are correlated with management and commercial efficiency - an efficient campaign impacts the final $\gamma$ - and correlated together: the higher the incentive $\delta$, the higher the probability of acceptance $\gamma$.\\
By fixing all other parameters and trying various combinations of $\delta$ and $\gamma$ we obtain the following 3D surfaces.\\
\begin{minipage}{\linewidth}
\makebox[\linewidth]{

  \includegraphics[width=0.7\linewidth]{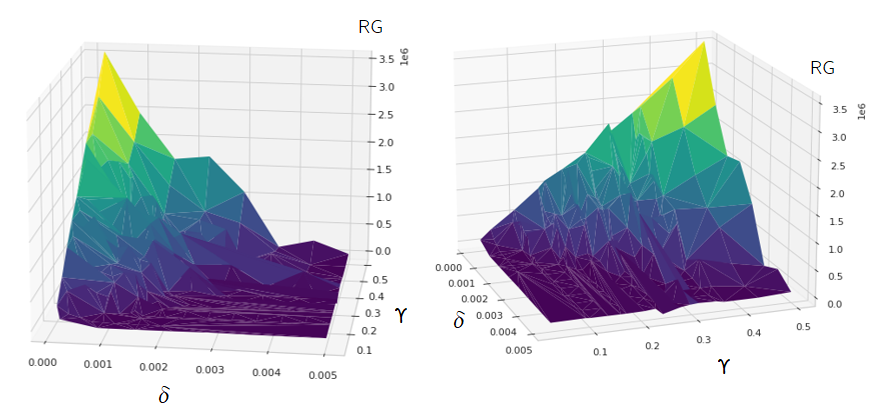}}
\captionof{figure}{3d plot ($\delta$, $\gamma$, RG)}\label{3d_plot_delta_gamma}
\end{minipage}\par\vspace{5mm}

This surface is not surprising and indicates that the higher the acceptance rate and the lower the incentive, the higher the retention gain. The surface gradient can give powerful insights regarding the most efficient commercial efforts to make: is it better for the insurer to propose lower incentives and manage to conserve the same acceptance probability or to put commercial effort into improving the acceptance probability for the same proposed incentive? This surface directly addresses this question.
\begin{mdframed}[hidealllines=true]
\begin{remark}
Of course, the interdependency of those parameters should make some part of this surface unrealistic from a management decision-making point of view. The insurer should consider such dependencies when designing a lapse management strategy.
\end{remark}
\end{mdframed}
Among the six parameters are also three conjectural parameters that depend on the external state of the world: the insurer's profitability $p$ (that depends on competition, macroeconomic considerations or regulation), the discount rate $d$ and the time horizon $T$ (that can be driven by the insurer's vision but also by regulation: the ORSA time horizon with the strategic and the long-term business planning time horizon should be both considered). Among them, we chose to fix $p$ and let $d$ and $T$ vary. Moreover, $p$ and $T$ are obviously interdependent and considered through the management's prospective view of the conjecture's evolution. A given interest rate scenario should represent a curve on the following surface.\\
\begin{minipage}{\linewidth}
\makebox[\linewidth]{

  \includegraphics[width=0.7\linewidth]{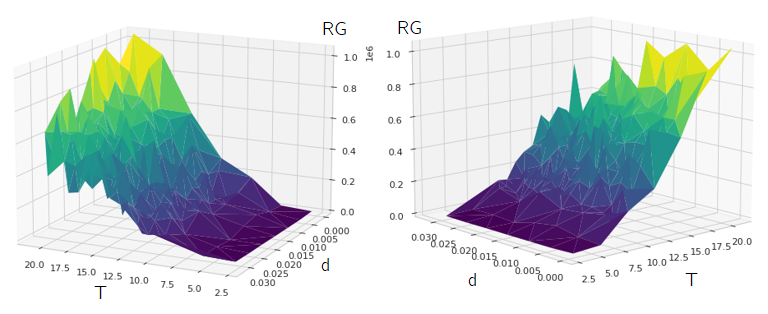}}
\captionof{figure}{3d plot (d, T, RG)}\label{3d_plot_d_T}
\end{minipage}
This surface is less smooth than the one displayed in Figure~\ref{3d_plot_delta_gamma} and seems to indicate a more unstable relationship between RG and the conjectural parameters. An explanation of that behavior can be that those surface points are generated by running our framework on a random subsample of our dataset, for computation time considerations. Generating the same surface with more policyholders is likely to give a smoother behavior.
\begin{mdframed}[hidealllines=true]
\begin{remark}
Of course, the interdependency of $T$ and $d$ should make some part of this surface unrealistic from an actuarial point of view. Actuarial rate projections would give precise plausible scenarios on this surface. Such considerations should be taken into account by the insurer when designing a lapse management strategy.
\end{remark}

\begin{remark}
The insurer can also use our framework to measure the retention gain to be expected at different time horizons obtained by existing retention campaigns. In that case, the insurer would have to neutralize the effect of the existing LMS in order to estimate the control portfolio's future value. We leave this remark as future work for applied risk management research.
\end{remark}
\end{mdframed}

\section{Conclusion and perspectives}\label{Conclusion}
The work carried out in this paper shows that including CLV in lapse management strategy can largely benefit an insurer's decision-making ability regarding lapse management strategy. We showed that survival tree-based models can outperform parametric approaches in such actuarial contexts. Then, our comparison of tree-based models on different lapse management strategies indicated that our CLV-based framework leads to increased predicted gains for any realistic scenario and acts as a loss-limiting targeting approach, regardless of the retention strategy. Moreover, the global results obtained in Section~\ref{General} show that our approach significantly improves on existing ones. Eventually, the discussion section highlighted the fact that our model can give insights to the life insurer regarding commercial and strategic decision-making.\\
\\
The framework and methodologies described in this paper suffer some limitations. For instance, following one single fixed strategy for every policyholder is arguably unrealistic. We could imagine an extension of our models to individualized lapse management strategies that would vary between subjects and could also be adjusted with time. In the application, we also considered constant parameters $p$ and $\delta$: a limiting assumption whose impact could be studied. There is also room for improvement regarding the correlations of LMS parameters: the value of the incentive and the acceptance probability are evidently interdependent parameters for an insurance company, and this interdependency could be considered. \\
This paper defines a practical management tool for life insurers as those models can measure the $RG$ and improve real strategies used in existing retention campaigns. Finally, our vision of CLV, and by extension, our whole methodology design could be improved by using longitudinal data that would yield time-dynamic results. We leave those two last observations for future work.\\
A real-life comparison between an actual retention strategy targeting and both the naïve and CLV-improved methodologies could be insightful for the insurer.

\section*{Acknowledgments, Statements and Declarations}
Work(s) conducted within the Research Chair DIALog under the aegis of the Risk Foundation, an initiative by CNP Assurances. The authors would like to express their very great gratitude to Marie Hyvernaud and Stéphanie Dosseh for their valuable and constructive suggestions while developing this research work. Special thanks should be given to Marie Hyvernaud for her contribution to code writing.

\newpage

{\small
\bibliography{bibfile}}

\newpage
\appendix
\section{Appendix}

\subsection{Competing risk framework}\label{appendix_competing_risks}
There are several regression models to estimate the global hazard and the hazard of one risk in settings where competing risks are present: modeling the cause-specific hazard and the subdistribution hazard function. They account for competing risks differently, obtaining different hazard functions and thus distinct advantages, drawbacks, and interpretations. Here, we will introduce those approaches' theoretical and practical implications and justify which one we will use in our modeling approaches.\\

In Cause-specific regression, each cause-specific hazard is estimated separately, in our case, the cause-specific hazards of lapse and death, by considering all subjects that experienced the competing event as censored. Here, $t$ is the traditional time variable of a survival model, with $t=0$ being the beginning of a policy. It is not to be confused with the use of $t$ in Sections~\ref{sec:Framework} and \ref{sec:Methodology}.
We remind that $J_{T}=0$ corresponds to an active subject that did not experience lapse $J_{T}=1$ or death $J_{T}=2$. The cause-specific hazard rates regarding the $j$-th risk ($j \in [1, \dots J]$) are defined as
$$
\lambda_{T, j}(t)=\lim _{d t \rightarrow 0} \frac{P\left(t \leq T<t+d t, J_{T}=j \mid T \geq t\right)}{d t} .
$$
We can recover the global hazard rate as $ \lambda_{T, 1}(t)+\cdots+\lambda_{T, J}(t)=\lambda_{T}(t)$, and derive the global survival distribution of $T$ as

\begin{align*}
P(T>t)&=1-F_{T}(t)=S_{T}(t)\\
&=\exp \left(-\int_{0}^{t}\left(\lambda_{T, 1}(s)+\cdots+\lambda_{T, J}(s)\right) d s\right) .
\end{align*}

This approach aims at analysing the cause-specific ``distribution" function: $F_{T, j}(t)=P\left(T \leq t, J_{T}=j\right)$. In practice, it is called the Cumulative Incidence Function ($CIF$) for cause $j$ and not a distribution function since $F_{T, j}(t) \rightarrow P\left(J_{T}=j\right) \neq 1 $ as $t \rightarrow+\infty$. By analogy with the classical survival framework, the $CIF$ can be characterised as $F_{T, j}(t)=\int_{0}^{t} f_{T, j}(s) d s$\footnote{We suppose that $T$ has a continuous distribution}, where $f_{T, j}$ is the improper\footnote{Because derived from the $CIF$, an improper cumulative distribution function} density function for cause $j$. It follows that
$$
f_{T, j}(s)=\lim _{d t \rightarrow 0} \frac{P\left(t \leq T<t+d t, J_{T}=j\right)}{d t}=\lambda_{T, j}(t) S_{T}(t) .
$$
The equation above is self-explanatory: the probability of experiencing cause $j$ at time $t$ is simply the product of surviving the previous time periods by the cause-specific hazard at time $t$. We finally obtain the $CIF$ for cause $j$ as
$$
F_{T, j}(t)=\int_{0}^{t} \lambda_{T, j}(s) \exp \left(-\int_{0}^{s} \lambda_{T}(u) d u\right) d s .
$$
There are several advantages to that approach. First of all, cause-specific hazard models can be easily fit with any classical implementation of CPH by simply considering as censored any subject that experienced the competing event. Then the $CIF$ is clearly interpretable and summable $P(T \leq t)=F_{T, 1}(s)+\cdots+F_{T, J}(s)$\footnote{unlike to the function $1-\exp \left(-\int_{0}^{t} \lambda_{T, j}(u) d u\right)$.}. On the other hand, the $CIF$ estimation of one given cause depends on all other causes: it implies that the study of a specific cause requires estimating the global hazard rate, and interpreting the effects of covariates on this cause is difficult. Indeed, part of the effects on a specific cause comes from the competing causes, but in our setting, we are only interested in the prediction of the survival probabilities, not their interpretation as such.  \\
\\
We have introduced it at  the beginning of this section; another approach is often considered to analyze competing risks and derive a cause-specific $CIF$. This other approach called the subdistribution hazard function of Fine and Gray regression, works by considering a new competing risk process $\tau$. Without loss of generality, let's consider death as our cause of interest, 
$$
\tau=T \times \mathbb{1}_{J_{T}=2}+\infty \times \mathbb{1}_{J_{T} \neq 2} .
$$
It has the same as $T$ regarding the risk of death, $P(\tau \leq t)=F_{T, 2}(t)$ and a mass point at infinity $1-F_{T, 2}(\infty)$, probability to observe other causes $\left(J_{T} \neq 2\right)$ or not to observe any failure. In other words, if the previous approach considered every subject that experienced competing events as censored, this approach considers a new and
artificial at-risk population. This last consideration is made clear when deriving the hazard rate of $\tau$,\\
$$
\lambda_{\tau}(t)=\lim _{d t \rightarrow 0} \frac{P\left(t \leq T<t+d t, J_{T}=2 \mid\{T \geq t\} \cup\left\{T \leq t, J_{T} \neq 2\right\}\right)}{d t} .
$$
Finally, we obtain the $CIF$ for the risk of death as
$$
F_{T, 2}(t)=1-\exp \left(-\int_{0}^{t} \lambda_{\tau}(s) d s\right) .
$$
This subdistribution approach resolves the most important drawback to cause-specific regression, as the coefficients resulting from it do have a direct relationship with the cumulative incidence: estimating the $CIF$ for a specific cause does not depend on the other causes, which makes the interpretation of $CIF$ easier. The subdistribution hazard models can be fit in R by using the crr function in the cmprsk package or using the timereg package. Still, to our knowledge, there is no implementation of a Fine and Gray model in Lifelines or, more generally, Python. We can also note that these two approaches are linked,\cite{putter2020a} and the link between $\lambda_{\tau}(t)$ and $\lambda_{T, j}(t)$ is given by
$$
\lambda_{\tau}(t)=r_{j}(t) \lambda_{T, j}(t), \text { with } r_{j}(t)=\frac{P(J_{T}=0)}{\displaystyle \sum^{J}_{p \neq j} P(J_{T}=p)} .
$$
In other words, if the probability of any competing risk is low, the two approaches give very close results.\\

\subsection{Survival analysis results}\label{appendix_surv}
The quantity $r^{lapser}_{i,t}$ represents the probability that the policy of subject $i$ is still active at time $t$, given that it was active at its last observed time. Predicting the overall conditional survival with the competing risks, in that case, can be achieved by creating a combined outcome. The policy ends with death or lapse, whichever comes first, and to compute $r^{lapser}$, we recode the competing events as a combined event. In terms of statistical guarantees, this approach is compatible with any survival analysis method.\\
In the following sections of this appendix, $r^{acceptant}_{i,t}$ indicates the probability of survival for subject $i$ at time $t$ given that it will not lapse. In other words, it is the survival probability regarding only the risk of death. As detailed in Section~\ref{com_risk_fram}, this corresponds to the cause-specific survival probability for death. It is to be noted that the density from which we derive our survival probabilities is improper as it derives itself from the $CIF$, which is not a proper distribution function.\footnote{as it does not tend to 1 as $t$ goes to $+\infty$}. Therefore, any conclusion about those probabilities should be drawn with care. Similarly to $r^{lapser}$, covariates selection and tuning are performed by minimizing AIC.\\
All graphs representing survival curves below are plotted with the same axis. The x-axes are the time in years, the y-axes represent the survival probability.

\subsubsection{Cox-model}
We first decide to estimate survival with a Cox Proportional hazard model with a spline baseline hazard from the Python library Lifelines. Covariate selection and tuning are performed by minimizing AIC. 
Here is what $r^{acceptant}$, the vector of cause-specific probabilities, looks like, and we can compare it to $r^{lapser}$ on some subjects.\\
\begin{minipage}{0.5\linewidth}
\makebox[\linewidth]{

  \includegraphics[width=0.95\linewidth]{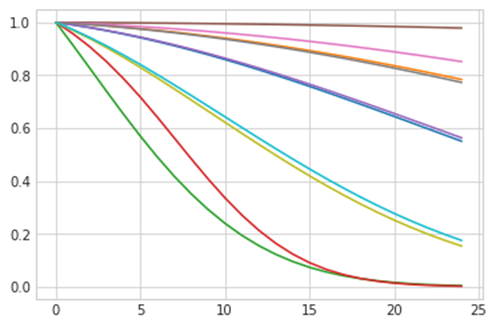}}
        \captionof{figure}{10 policyholders' survival curve for $r^{acceptant}$ with Cox model}
        \label{fig:plot_surv_10_r_stay}
\end{minipage}
\begin{minipage}{0.5\linewidth}
\makebox[\linewidth]{
  \includegraphics[width=0.95\linewidth]{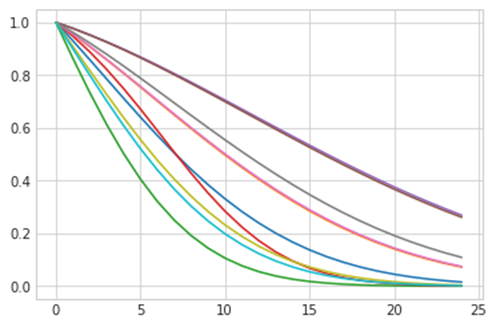}}
        \captionof{figure}{10 policyholders' survival curve for $r^{lapser}$}
        \label{fig:fig:plot_surv_10_r_lapse}
\end{minipage}

The effect of various covariates on the survival outcome can be found below.\\
\begin{minipage}{\linewidth}
\makebox[\linewidth]{
  \includegraphics[width=0.45\linewidth]{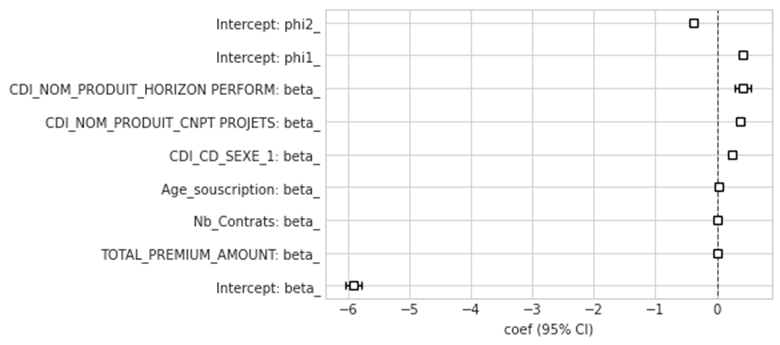}}
        \captionof{figure}{Coefficient plot for $r^{lapser}$}
        \label{fig:plot_coef_r_lapse}
\end{minipage}
\vspace{10mm}\\
\begin{minipage}{.5\linewidth}
\makebox[\linewidth]{
  \includegraphics[width=0.95\linewidth]{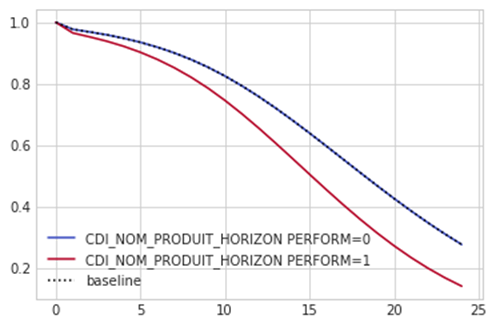}}
        \captionof{figure}{$r^{lapser}$ trajectories for different products}
        \label{fig:plot_surv_product_r_lapse}
\end{minipage}
\begin{minipage}{.5\linewidth}
\makebox[\linewidth]{
  \includegraphics[width=0.95\linewidth]{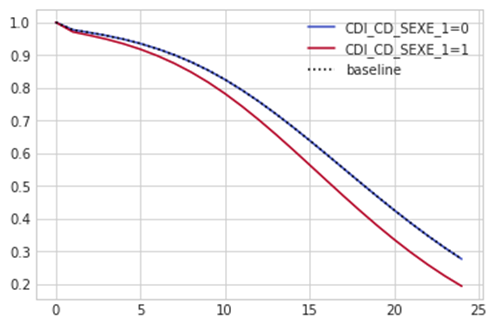}}
        \captionof{figure}{$r^{lapser}$ trajectories by gender}
        \label{fig:plot_surv_sex_r_lapse}
\end{minipage}
\begin{minipage}{.5\linewidth}
\makebox[\linewidth]{
  \includegraphics[width=0.95\linewidth]{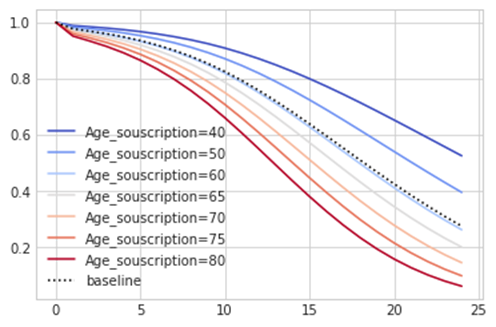}}
        \captionof{figure}{$r^{lapser}$ trajectories  for different ages}
        \label{fig:plot_surv_age_r_lapse}
\end{minipage}
\begin{minipage}{.5\linewidth}
\makebox[\linewidth]{
  \includegraphics[width=0.95\linewidth]{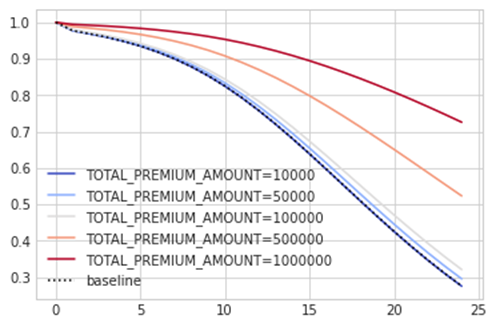}}
        \captionof{figure}{$r^{lapser}$ trajectories  for different face amounts}
        \label{fig:plot_surv_outstanding_amount_r_lapse}
\end{minipage}
\phantom{empty}
\begin{minipage}{.5\linewidth}
\makebox[\linewidth]{
  \includegraphics[width=0.95\linewidth]{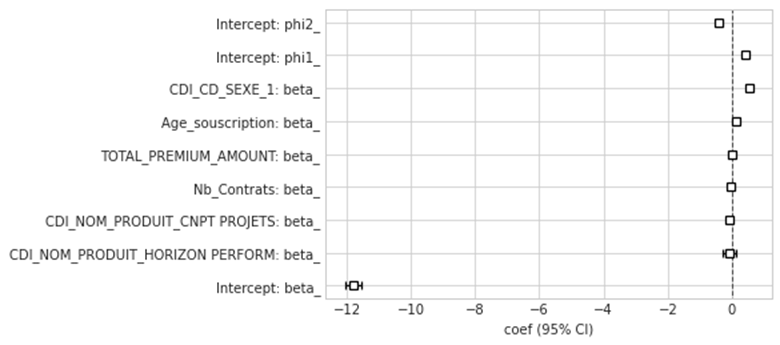}}
        \captionof{figure}{Coefficient plot for $r^{acceptant}$}
        \label{fig:plot_coef_r_stay}
\end{minipage}
\begin{minipage}{.5\linewidth}
\makebox[\linewidth]{
  \includegraphics[width=0.95\linewidth]{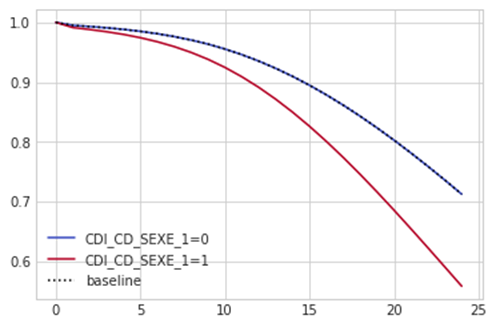}}
        \captionof{figure}{$r^{acceptant}$ trajectories by gender}
        \label{fig:plot_surv_sex_r_stay}
\end{minipage}

\phantom{empty}
\begin{minipage}{.5\linewidth}
\makebox[\linewidth]{
  \includegraphics[width=0.95\linewidth]{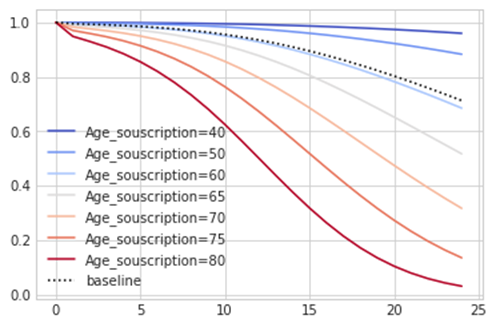}}
\captionof{figure}{$r^{acceptant}$ trajectories for different ages}
        \label{fig:plot_surv_age_r_stay}
\end{minipage}
\begin{minipage}{.5\linewidth}
\makebox[\linewidth]{
  \includegraphics[width=0.95\linewidth]{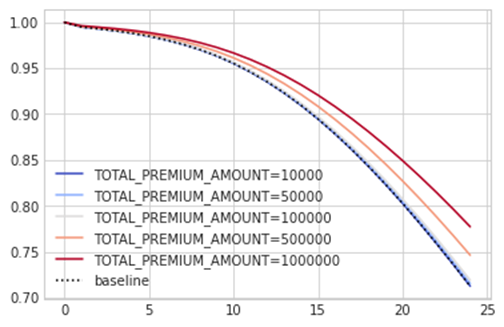}}
\captionof{figure}{$r^{acceptant}$ trajectories for different face amounts}\label{fig:plot_surv_outstanding_amount_r_stay}
\end{minipage}

\subsubsection{RSF}\label{appendix_rsf}
We obtain better results than Cox in terms of concordance index at the cost of very high computation time for one training with one set of parameters - 5days without parallelisation, 4 hours with - compared to a few seconds for cox model.\\
Some of the results we obtain are displayed below. \\

\begin{minipage}{.5\linewidth}
\makebox[\linewidth]{
  \includegraphics[width=0.95\linewidth]{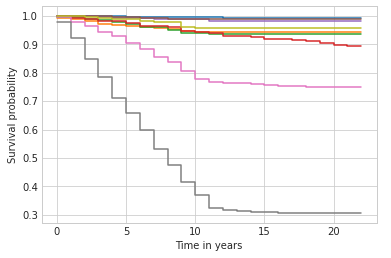}}
        \captionof{figure}{10 policyholders' survival curve for $r^{acceptant}$ with RSF}
        \label{fig:plot_surv_10_r_stay_rsf}
\end{minipage}
\begin{minipage}{.5\linewidth}
\makebox[\linewidth]{
  \includegraphics[width=0.95\linewidth]{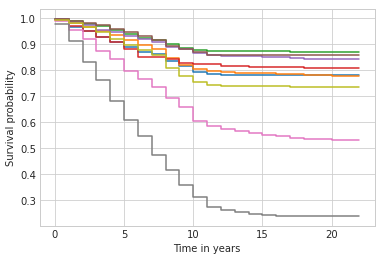}}
        \captionof{figure}{10 policyholders' survival curve for $r^{lapser}$ with RSF}
        \label{fig:plot_surv_10_r_lapse_rsf}
\end{minipage}

    \begin{tabular}{llll}
    \multicolumn{3}{c}{Weight} & \multicolumn{1}{c}{Feature} \\
    \midrule
    \rowcolor[rgb]{ 1,  1,  0} 0.3148 & ±     & 0.0064 & Age\_souscription \\
    \rowcolor[rgb]{ 1,  1,  0} 0.0100 & ±     & 0.0008 & CDI\_CD\_SEXE\_1 \\
    0.0091 & ±     & 0.0014 & PRODUIT\_2 \\
    0.0077 & ±     & 0.0006 & TOTAL\_PREMIUM\_AMOUNT \\
    0.0013 & ±     & 0.0004 & Nb\_Contrats \\
    0.0010 & ±     & 0.0003 & PRODUIT\_3 \\
    \end{tabular}%
    \captionof{table}{Covariates importance for $r^{acceptant}$ with RSF}\label{fig:variables_importance_rsf_r_stay}
    
    \begin{tabular}{llll}
    \multicolumn{3}{c}{Weight} & \multicolumn{1}{c}{Feature} \\
    \midrule
    \rowcolor[rgb]{ 1,  1,  0} 0.1838 & ±     & 0.0045 & Age\_souscription \\
    \rowcolor[rgb]{ 1,  1,  0} 0.0415 & ±     & 0.0018 & TOTAL\_PREMIUM\_AMOUNT \\
    0.0083 & ±     & 0.0011 & CDI\_CD\_SEXE\_1 \\
    0.0026 & ±     & 0.0013 & PRODUIT\_2 \\
    0.0022 & ±     & 0.0006 & PRODUIT\_3 \\
    0.0020 & ±     & 0.0006 & Nb\_Contrats \\
    \end{tabular}%
    \captionof{table}{Covariates importance for $r^{lapser}$ with RSF}\label{fig:variables_importance_rsf_r_lapse}

\subsubsection{XGSB}\label{appendix_gbsm}
We obtain better results than Cox and slightly better results than RSF in terms of concordance index at the cost of even higher computation time for one training with one set of parameters - 10h with great parallelisation - compared to a few seconds for Cox model.\\
Some of the results we obtain are displayed below. \\

\begin{minipage}{.5\linewidth}
\makebox[\linewidth]{

  \includegraphics[width=0.93\linewidth]{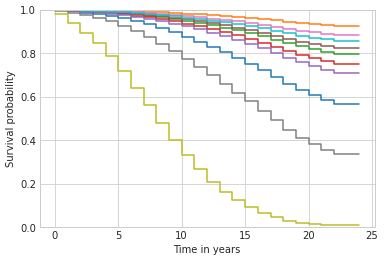}}
        \captionof{figure}{10 policyholders' survival curve for $r^{acceptant}$ with GBSM}
        \label{fig:plot_surv_10_r_stay_gbsm}
\end{minipage}
\begin{minipage}{.5\linewidth}
\makebox[\linewidth]{
  \includegraphics[width=0.93\linewidth]{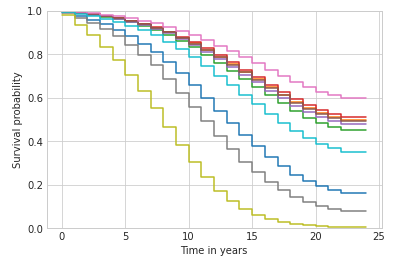}}
        \captionof{figure}{10 policyholders' survival curve for $r^{lapser}$ with GBSM}
        \label{fig:plot_surv_10_r_lapse_gbsm}
\end{minipage}

    \begin{tabular}{llll}
    \multicolumn{3}{c}{Weight} & \multicolumn{1}{c}{Feature} \\
    \midrule
    \rowcolor[rgb]{ 1,  1,  0} 0.3274 & ±     & 0.0071 & Age\_souscription \\
    \rowcolor[rgb]{ 1,  1,  0} 0.0104 & ±     & 0.0006 & TOTAL\_PREMIUM\_AMOUNT \\
    \rowcolor[rgb]{ 1,  1,  0} 0.0100 & ±     & 0.0008 & CDI\_CD\_SEXE\_1 \\
    0.0025 & ±     & 0.0005 & PRODUIT\_2 \\
    0.0005 & ±     & 0.0001 & Nb\_Contrats \\
    0.0000 & ±     & 0.0001 & PRODUIT\_3 \\
    \end{tabular}%
    \captionof{table}{Covariates importance for $r^{acceptant}$ with GBSM}\label{fig:variables_importance_gbsm_r_stay}

    \begin{tabular}{llll}
    \multicolumn{3}{c}{Weight} & \multicolumn{1}{c}{Feature} \\
    \midrule
    \rowcolor[rgb]{ 1,  1,  0} 0.1872 & ±     & 0.0039 & Age\_souscription \\
    \rowcolor[rgb]{ 1,  1,  0} 0.0438 & ±     & 0.0020 & TOTAL\_PREMIUM\_AMOUNT \\
    \rowcolor[rgb]{ 1,  1,  0} 0.0134 & ±     & 0.0014 & PRODUIT\_2 \\
    0.0076 & ±     & 0.0009 & CDI\_CD\_SEXE\_1 \\
    0.0051 & ±     & 0.0006 & PRODUIT\_3 \\
    0.0011 & ±     & 0.0004 & Nb\_Contrats \\
    \end{tabular}%
    \captionof{table}{Covariates importance for $r^{lapser}$ with GBSM}\label{fig:variables_importance_gbsm_r_lapse}    
    
\subsubsection{Final survival model}
The final concordance index scores are displayed below:\\

\begin{minipage}{\linewidth}
\makebox[\linewidth]{
  \begin{tabular}{r|cc}
    \multicolumn{1}{r}{} & \multicolumn{2}{c}{Concordance Index} \\
\cmidrule{2-3}    \multicolumn{1}{r}{} & $r^{lapser}$ & $r^{acceptant}$ \\
\cmidrule{2-3}    Cox model & 69,5\% & 80,7\% \\
    RSF   & 71,6\% & 83,7\% \\
    GBSM  & 73,0\% & 84,1\% \\
    \end{tabular}%
    }
\captionof{table}{Survival models comparison}\label{Survival_models_comparison}
\end{minipage}

\subsection{Other results}\label{other_results}
\begin{minipage}{1\linewidth}
\makebox[\linewidth]{

  \includegraphics[width=0.75\linewidth]{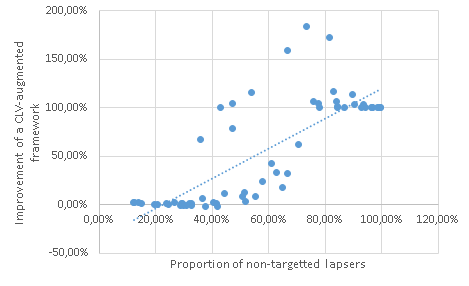}
  }
        \captionof{figure}{Correlation between the proportion of non-targeted lapsers and the improvement\protect\footnotemark of a CLV-augmented LMS}
        \label{fig:improvement_correlation}
\end{minipage}

\subsection{Considering various statistical metrics}\label{other_metrics}

The table below contains the results of the " LMS listed in Table~\ref{tab:Numerical_results_1}, evaluated on accuracy, recall, F1-score and AUC. For every metric, it displays the results of a classification over $y_i$ tuned and cross-validated with each of the metrics - respectively $\stackrel{accuracy}{y_i}$, $\stackrel{recall}{y_i}$, $\stackrel{F1-score}{y_i}$ and $\stackrel{AUC}{y_i}$ - or over $\tilde{y}_i$ which is always tuned and cross-validated with $RG$.\\
\\
\begin{adjustbox}{width=1\textwidth}
    \begin{tabular}{ccrrrrrrrrrrrrrrrrrr}
        \toprule
        \multirow{2}{*}{N°} & \multirow{2}{*}{Model} & \multicolumn{2}{c}{Accuracy} & 
        \multicolumn{2}{c}{Recall} & 
        \multicolumn{2}{c}{F1-score} & 
        \multicolumn{2}{c}{AUC} & \multicolumn{5}{c}{Retention gain} & \multicolumn{5}{c}{RG/target} \\
        \cmidrule(lr){3-4} \cmidrule(lr){5-6} \cmidrule(lr){7-8} \cmidrule(lr){9-10} \cmidrule(lr){11-15} \cmidrule(lr){16-20}
        & & \multicolumn{1}{c}{$\stackrel{accuracy}{y_i}$} & \multicolumn{1}{c}{$\tilde{y}_i$} & \multicolumn{1}{c}{$\stackrel{recall}{y_i}$} & \multicolumn{1}{c}{$\tilde{y}_i$} & \multicolumn{1}{c}{$\stackrel{F1-score}{y_i}$} & \multicolumn{1}{c}{$\tilde{y}_i$} & \multicolumn{1}{c}{$\stackrel{AUC}{y_i}$} & \multicolumn{1}{c}{$\tilde{y}_i$} & \multicolumn{1}{c}{$\stackrel{accuracy}{y_i}$} & \multicolumn{1}{c}{$\stackrel{recall}{y_i}$} & \multicolumn{1}{c}{$\stackrel{F1-score}{y_i}$} & \multicolumn{1}{c}{$\stackrel{AUC}{y_i}$} & \multicolumn{1}{c}{$\tilde{y}_i$} & \multicolumn{1}{c}{$\stackrel{accuracy}{y_i}$} & \multicolumn{1}{c}{$\stackrel{recall}{y_i}$} & \multicolumn{1}{c}{$\stackrel{F1-score}{y_i}$} & \multicolumn{1}{c}{$\stackrel{AUC}{y_i}$} & \multicolumn{1}{c}{$\tilde{y}_i$} \\
        \midrule
        \multirow{3}{*}{A-1} & CART  & 92.3\% & 85.3\% & 75.7\% & 18.0\% & 76.4\% & 28.9\% & 85.6\% & 58.4\% & 114,661 & 128,251 & 135,085 & 128,251 & 219,655 & 4.48 & 5.13 & 5.46 & 5.13 & 38.20 \\
         & RF & 92.9\% & 85.4\% & 74.7\% & 14.7\% & 77.7\% & 24.9\% & 85.6\% & 57.0\% & 232,314 & 203,821 & 203,821 & 203,821 & 287,884 & 9.82 & 8.66 & 8.66 & 8.66 & 56.65 \\
         & XGB & 93.4\% & 85.8\% & 79.6\% & 18.5\% & 80.0\% & 30.1\% & 87.8\% & 58.8\% & 243,365 & 261,553 & 266,198 & 266,198 & 324,952 & 9.61 & 10.25 & 10.55 & 10.55 & 54.64 \\
        \midrule
        \multirow{3}{*}{A-5} & CART  & 92.3\% & 83.6\% & 76.8\% & 2.5\% & 76.7\% & 4.8\% & 86.0\% & 51.1\% & -514,477  & -497,277 & -494,414 & -497,277 &  - 112,372 & - 20.08 & -19.28 & -19.21 & -19.28 & - 86.48 \\
        & RF & 92.9\% & 83.4\% & 74.9\% & 0.6\% & 77.9\% & 1.3\% & 85.7\% & 50.3\% & - 323,544 & -323,543 & -325,110 & -323,543 & -3,937 & - 13.65 & - 13.64 & - 13.72 & - 13.64 & -28.28 \\
        & XGB & 93.4\% & 83.3\% & 79.8\% & 0.4\% & 80.0\% & 0.7\% & 87.9\% & 50.2\% & -  383,004 & -379,736 & -379,736 & -379,736 & 0 & - 15.14 & -14.88 & -14.88 & -14.88 & 0 \\
        \midrule
        \multirow{3}{*}{A-25} & CART  & 92.3\% & 89.2\% & 76.2\% & 50.2\% & 76.7\% & 60.1\% & 85.8\% & 73.6\% & 4,160,423 & 4,322,030 & 4,267,310 & 4,322,030 & 3,882,623 & 162.44 & 170.53 & 251.34 & 251.34 & 241.06 \\
        & RF & 92.9\% & 89.5\% & 75.3\% & 47.5\% & 78.0\% & 60.3\% & 85.8\% & 72.7\% & 4,018,432 & 3,687,705 & 3,687,705 & 3,687,705 & 3,666,219 & 169.65 & 154.49 & 154.49 & 154.49 & 249.54 \\
        & XGB   & 93.4\% &90.0\% & 80.1\% & 51.2\% & 80.1\% & 63.0\% & 88.0\% & 74.4\% & 4,455,108 & 4,633,404 & 4,578,684 & 4,633,404 & 4,410,629 & 176.09 & 179.36 & 180.16 & 179.36 & 267.87 \\
        \bottomrule
    \end{tabular}
\end{adjustbox}   
\captionof{table}{Results of representative LMS with various statistical metrics}\label{other_metrics_results}%
It is to be noted that regardless of the evaluation metric used for tuning and validation purposes, the objective function used with XGB to generate those results is always the log-loss function. Using the area under the ROC curve or the area under the Precision-Recall curve as an objective function in this boosting algorithm would surely yield better results when trained on $y_i$ and even better on the more unbalanced $\tilde{y}_i$. As stated in Section~\ref{classif_step}, this analysis is not within the scope of our article.

\footnotetext{Taking the results of XGBoost and excluding LMS n°B-27 that has a very high improvement ratio.}

\subsection{Complete LMS numerical results}\label{more_lms_results}

\begin{adjustbox}{width=1\textwidth}
\begin{minipage}{.5\linewidth}

\begin{adjustbox}{width=0.85\textwidth}
    \begin{tabular}{c|cccccc}
    \multicolumn{1}{c}{LMS} & p     & $\delta$ & $\gamma$ & c     & d     & T \\
    \midrule
    A-1   & 2,50\% & 0,04\% & 25\%  & 10    & 1,50\% &                     5  \\
    A-2   & 2,50\% & 0,04\% & 25\%  & 10    & 1,50\% &                   20  \\
    A-3   & 2,50\% & 0,04\% & 25\%  & 100   & 1,50\% &                     5  \\
    A-4   & 2,50\% & 0,04\% & 25\%  & 100   & 1,50\% &                   20  \\
    A-5   & 2,50\% & 0,04\% & 5\%   & 10    & 1,50\% &                     5  \\
    A-6   & 2,50\% & 0,04\% & 5\%   & 10    & 1,50\% &                   20  \\
    A-7   & 2,50\% & 0,04\% & 5\%   & 100   & 1,50\% &                     5  \\
    A-8   & 2,50\% & 0,04\% & 5\%   & 100   & 1,50\% &                   20  \\
    A-9   & 2,50\% & 0,10\% & 25\%  & 10    & 1,50\% &                     5  \\
    A-10  & 2,50\% & 0,10\% & 25\%  & 10    & 1,50\% &                   20  \\
    A-11  & 2,50\% & 0,10\% & 25\%  & 100   & 1,50\% &                     5  \\
    A-12  & 2,50\% & 0,10\% & 25\%  & 100   & 1,50\% &                   20  \\
    A-13  & 2,50\% & 0,10\% & 5\%   & 10    & 1,50\% &                     5  \\
    A-14  & 2,50\% & 0,10\% & 5\%   & 10    & 1,50\% &                   20  \\
    A-15  & 2,50\% & 0,10\% & 5\%   & 100   & 1,50\% &                     5  \\
    A-16  & 2,50\% & 0,10\% & 5\%   & 100   & 1,50\% &                   20  \\
    A-17  & 5,00\% & 0,04\% & 25\%  & 10    & 1,50\% &                     5  \\
    A-18  & 5,00\% & 0,04\% & 25\%  & 10    & 1,50\% &                   20  \\
    A-19  & 5,00\% & 0,04\% & 25\%  & 100   & 1,50\% &                     5  \\
    A-20  & 5,00\% & 0,04\% & 25\%  & 100   & 1,50\% &                   20  \\
    A-21  & 5,00\% & 0,04\% & 5\%   & 10    & 1,50\% &                     5  \\
    A-22  & 5,00\% & 0,04\% & 5\%   & 10    & 1,50\% &                   20  \\
    A-23  & 5,00\% & 0,04\% & 5\%   & 100   & 1,50\% &                     5  \\
    A-24  & 5,00\% & 0,04\% & 5\%   & 100   & 1,50\% &                   20  \\
    A-25  & 5,00\% & 0,10\% & 25\%  & 10    & 1,50\% &                     5  \\
    A-26  & 5,00\% & 0,10\% & 25\%  & 10    & 1,50\% &                   20  \\
    A-27  & 5,00\% & 0,10\% & 25\%  & 100   & 1,50\% &                     5  \\
    A-28  & 5,00\% & 0,10\% & 25\%  & 100   & 1,50\% &                   20  \\
    A-29  & 5,00\% & 0,10\% & 5\%   & 10    & 1,50\% &                     5  \\
    A-30  & 5,00\% & 0,10\% & 5\%   & 10    & 1,50\% &                   20  \\
    A-31  & 5,00\% & 0,10\% & 5\%   & 100   & 1,50\% &                     5  \\
    A-32  & 5,00\% & 0,10\% & 5\%   & 100   & 1,50\% &                   20  \\
    \end{tabular}%
    \end{adjustbox}
\end{minipage}
\begin{minipage}{.5\linewidth}
\begin{adjustbox}{width=0.85\textwidth}
    \begin{tabular}{c|cccccc}
    \multicolumn{1}{c}{LMS} & p     & $\delta$ & $\gamma$ & c     & d     & T \\
    \midrule
    B-1   & 2,50\% & 0,08\% & 20\%  & 10    & 1,50\% &                     5  \\
    B-2   & 2,50\% & 0,08\% & 20\%  & 10    & 1,50\% &                   20  \\
    B-3   & 2,50\% & 0,08\% & 20\%  & 100   & 1,50\% &                     5  \\
    B-4   & 2,50\% & 0,08\% & 20\%  & 100   & 1,50\% &                   20  \\
    B-5   & 2,50\% & 0,08\% & 10\%  & 10    & 1,50\% &                     5  \\
    B-6   & 2,50\% & 0,08\% & 10\%  & 10    & 1,50\% &                   20  \\
    B-7   & 2,50\% & 0,08\% & 10\%  & 100   & 1,50\% &                     5  \\
    B-8   & 2,50\% & 0,08\% & 10\%  & 100   & 1,50\% &                   20  \\
    B-9   & 2,50\% & 0,20\% & 20\%  & 10    & 1,50\% &                     5  \\
    B-10  & 2,50\% & 0,20\% & 20\%  & 10    & 1,50\% &                   20  \\
    B-11  & 2,50\% & 0,20\% & 20\%  & 100   & 1,50\% &                     5  \\
    B-12  & 2,50\% & 0,20\% & 20\%  & 100   & 1,50\% &                   20  \\
    B-13  & 2,50\% & 0,20\% & 10\%  & 10    & 1,50\% &                     5  \\
    B-14  & 2,50\% & 0,20\% & 10\%  & 10    & 1,50\% &                   20  \\
    B-15  & 2,50\% & 0,20\% & 10\%  & 100   & 1,50\% &                     5  \\
    B-16  & 2,50\% & 0,20\% & 10\%  & 100   & 1,50\% &                   20  \\
    B-17  & 5,00\% & 0,08\% & 20\%  & 10    & 1,50\% &                     5  \\
    B-18  & 5,00\% & 0,08\% & 20\%  & 10    & 1,50\% &                   20  \\
    B-19  & 5,00\% & 0,08\% & 20\%  & 100   & 1,50\% &                     5  \\
    B-20  & 5,00\% & 0,08\% & 20\%  & 100   & 1,50\% &                   20  \\
    B-21  & 5,00\% & 0,08\% & 10\%  & 10    & 1,50\% &                     5  \\
    B-22  & 5,00\% & 0,08\% & 10\%  & 10    & 1,50\% &                   20  \\
    B-23  & 5,00\% & 0,08\% & 10\%  & 100   & 1,50\% &                     5  \\
    B-24  & 5,00\% & 0,08\% & 10\%  & 100   & 1,50\% &                   20  \\
    B-25  & 5,00\% & 0,20\% & 20\%  & 10    & 1,50\% &                     5  \\
    B-26  & 5,00\% & 0,20\% & 20\%  & 10    & 1,50\% &                   20  \\
    B-27  & 5,00\% & 0,20\% & 20\%  & 100   & 1,50\% &                     5  \\
    B-28  & 5,00\% & 0,20\% & 20\%  & 100   & 1,50\% &                   20  \\
    B-29  & 5,00\% & 0,20\% & 10\%  & 10    & 1,50\% &                     5  \\
    B-30  & 5,00\% & 0,20\% & 10\%  & 10    & 1,50\% &                   20  \\
    B-31  & 5,00\% & 0,20\% & 10\%  & 100   & 1,50\% &                     5  \\
    B-32  & 5,00\% & 0,20\% & 10\%  & 100   & 1,50\% &                   20  \\
    \end{tabular}%
    \end{adjustbox}

\end{minipage}
\end{adjustbox}
\captionof{table}{More LMS}\label{more_lms_considered2}

\begin{adjustbox}{width=1\textwidth}
    \begin{tabular}{cclcrrrrrrr}
    \multirow{2}[1]{*}{N°} & \multirow{2}[1]{*}{time (s)} & \multicolumn{1}{c}{\multirow{2}[1]{*}{Model}} & \multicolumn{1}{c}{\multirow{2}[1]{*}{\% target diff}} & \multicolumn{2}{c}{Accuracy} & \multicolumn{2}{c}{ Retention gain } & \multicolumn{2}{c}{ RG/target } & \multicolumn{1}{c}{\multirow{2}[1]{*}{Improvement\footref{improvement_note}}} \\
          &       &       &       & \multicolumn{1}{c}{$y_i$} & \multicolumn{1}{c}{$\tilde{y}_i$} & \multicolumn{1}{c}{ $y_i$ } & \multicolumn{1}{c}{ $\tilde{y}_i$ } & \multicolumn{1}{c}{ $y_i$ } & \multicolumn{1}{c}{ $\tilde{y}_i$ } &  \\
    \midrule
    \multirow{3}[2]{*}{A-1} & \multirow{3}[2]{*}{4949} & CART  & \multirow{3}[2]{*}{62,58\%} & 92,3\% & 85,3\% &                114 661  &                219 655  &                4,48  &              38,20  & 91,57\% \\
          &       & RF    &       & 92,9\% & 85,4\% &                232 314  &                287 884  &                9,82  &              56,65  & 23,92\% \\
          &       & XGB   &       & 93,4\% & 85,8\% &                243 365  &                324 952  &                9,61  &              54,64  & 33,52\% \\
    \midrule
    \multirow{3}[2]{*}{A-2} & \multirow{3}[2]{*}{6111} & CART  & \multirow{3}[2]{*}{26,66\%} & 92,3\% & 89,8\% &             7 092 097  &             6 142 119  &            277,00  &            353,83  & -13,39\% \\
          &       & RF    &       & 92,9\% & 90,2\% &             6 596 374  &             5 696 455  &            278,47  &            351,02  & -13,64\% \\
          &       & XGB   &       & 93,4\% & 90,9\% &             7 308 721  &             7 432 688  &            288,92  &            404,84  & 1,70\% \\
    \midrule
    \multirow{3}[2]{*}{A-3} & \multirow{3}[2]{*}{4603} & CART  & \multirow{3}[2]{*}{93,50\%} & 92,3\% & 83,3\% & -          2 187 622  & -                  8 224  & -           85,52  & -           31,09  & 99,62\% \\
          &       & RF    &       & 92,9\% & 83,4\% & -          1 900 265  &                  45 483  & -           80,18  &            194,35  & 102,39\% \\
          &       & XGB   &       & 93,4\% & 83,5\% & -          2 032 650  &                  77 481  & -           80,39  &            174,44  & 103,81\% \\
    \midrule
    \multirow{3}[2]{*}{A-4} & \multirow{3}[2]{*}{5555} & CART  & \multirow{3}[2]{*}{55,37\%} & 92,3\% & 86,5\% &             4 789 814  &             5 117 844  &            187,00  &            577,74  & 6,85\% \\
          &       & RF    &       & 92,9\% & 86,4\% &             4 463 796  &             4 255 175  &            188,47  &            566,05  & -4,67\% \\
          &       & XGB   &       & 93,4\% & 86,8\% &             5 032 706  &             5 433 366  &            198,92  &            610,26  & 7,96\% \\
    \midrule
    \multirow{3}[2]{*}{A-5} & \multirow{3}[2]{*}{4753} & CART  & \multirow{3}[2]{*}{86,72\%} & 92,3\% & 83,6\% & -             514 477  & -             112 372  & -           20,08  & -           86,48  & 78,16\% \\
          &       & RF    &       & 92,9\% & 83,4\% & -             323 544  & -                  3 937  & -           13,65  & -           28,28  & 98,78\% \\
          &       & XGB   &       & 93,4\% & 83,3\% & -             383 004  &                           0    & -           15,14  &                     0    & 100,00\% \\
    \midrule
    \multirow{3}[2]{*}{A-6} & \multirow{3}[2]{*}{5803} & CART  & \multirow{3}[2]{*}{44,27\%} & 92,3\% & 87,9\% &                335 810  &                517 224  &              13,17  &              39,91  & 54,02\% \\
          &       & RF    &       & 92,9\% & 87,9\% &                655 350  &                661 021  &              27,68  &              61,13  & 0,87\% \\
          &       & XGB   &       & 93,4\% & 88,6\% &                654 219  &                729 493  &              25,86  &              58,22  & 11,51\% \\
    \midrule
    \multirow{3}[2]{*}{A-7} & \multirow{3}[2]{*}{4241} & CART  & \multirow{3}[2]{*}{99,09\%} & 92,3\% & 83,3\% & -          2 816 759  & -                10 205  & -         110,08  & -         384,04  & 99,64\% \\
          &       & RF    &       & 92,9\% & 83,3\% & -          2 456 122  &                     1 013  & -         103,65  &              66,30  & 100,04\% \\
          &       & XGB   &       & 93,4\% & 83,3\% & -          2 659 020  &                        243  & -         105,14  &              15,92  & 100,01\% \\
    \midrule
    \multirow{3}[2]{*}{A-8} & \multirow{3}[2]{*}{5164} & CART  & \multirow{3}[2]{*}{82,78\%} & 92,3\% & 84,0\% & -          1 966 473  & -                46 323  & -           76,83  & -           22,31  & 97,64\% \\
          &       & RF    &       & 92,9\% & 84,0\% & -          1 477 229  &                253 885  & -           62,32  &            149,67  & 117,19\% \\
          &       & XGB   &       & 93,4\% & 84,1\% & -          1 621 796  &                273 243  & -           64,14  &            117,83  & 116,85\% \\
    \midrule
    \multirow{3}[2]{*}{A-9} & \multirow{3}[2]{*}{4781} & CART  & \multirow{3}[2]{*}{77,60\%} & 92,3\% & 83,7\% & -             825 372  & -             161 100  & -           32,19  & -         127,87  & 80,48\% \\
          &       & RF    &       & 92,9\% & 83,4\% & -             384 736  &                     8 596  & -           16,22  &              32,12  & 102,23\% \\
          &       & XGB   &       & 93,4\% & 83,6\% & -             498 263  &                  22 337  & -           19,70  &              35,47  & 104,48\% \\
    \midrule
    \multirow{3}[2]{*}{A-10} & \multirow{3}[2]{*}{6075} & CART  & \multirow{3}[2]{*}{29,10\%} & 92,3\% & 89,7\% &             4 614 513  &             4 483 831  &            180,36  &            266,33  & -2,83\% \\
          &       & RF    &       & 92,9\% & 89,9\% &             4 973 929  &             4 328 724  &            210,01  &            280,90  & -12,97\% \\
          &       & XGB   &       & 93,4\% & 90,7\% &             5 354 770  &             5 368 917  &            211,69  &            301,57  & 0,26\% \\
    \midrule
    \multirow{3}[2]{*}{A-11} & \multirow{3}[2]{*}{4506} & CART  & \multirow{3}[2]{*}{96,56\%} & 92,3\% & 83,2\% & -          3 127 655  & -             118 886  & -         122,19  & -      2 230,39  & 96,20\% \\
          &       & RF    &       & 92,9\% & 83,3\% & -          2 517 315  &                     1 340  & -         106,22  &              87,71  & 100,05\% \\
          &       & XGB   &       & 93,4\% & 83,3\% & -          2 774 278  &                        736  & -         109,70  &              52,00  & 100,03\% \\
    \midrule
    \multirow{3}[2]{*}{A-12} & \multirow{3}[2]{*}{5534} & CART  & \multirow{3}[2]{*}{57,93\%} & 92,3\% & 86,2\% &             2 312 231  &             3 310 314  &              90,36  &            412,71  & 43,17\% \\
          &       & RF    &       & 92,9\% & 86,1\% &             2 841 351  &             3 129 652  &            120,01  &            465,74  & 10,15\% \\
          &       & XGB   &       & 93,4\% & 86,6\% &             3 078 755  &             3 825 920  &            121,69  &            475,53  & 24,27\% \\
    \midrule
    \multirow{3}[2]{*}{A-13} & \multirow{3}[2]{*}{4640} & CART  & \multirow{3}[2]{*}{92,91\%} & 92,3\% & 83,3\% & -          1 201 626  & -             163 056  & -           46,87  & -      1 838,44  & 86,43\% \\
          &       & RF    &       & 92,9\% & 83,3\% & -             717 620  & -                  5 339  & -           30,28  & -         354,24  & 99,26\% \\
          &       & XGB   &       & 93,4\% & 83,3\% & -             875 378  &                        508  & -           34,60  &              16,26  & 100,06\% \\
    \midrule
    \multirow{3}[2]{*}{A-14} & \multirow{3}[2]{*}{5739} & CART  & \multirow{3}[2]{*}{47,12\%} & 92,3\% & 87,3\% & -          1 476 651  & -             831 019  & -           57,49  & -           77,99  & 43,72\% \\
          &       & RF    &       & 92,9\% & 86,0\% & -             380 683  &                126 532  & -           16,03  &              21,14  & 133,24\% \\
          &       & XGB   &       & 93,4\% & 85,5\% & -             644 389  &                  29 382  & -           25,47  &                7,10  & 104,56\% \\
    \midrule
    \multirow{3}[2]{*}{A-15} & \multirow{3}[2]{*}{4216} & CART  & \multirow{3}[2]{*}{99,61\%} & 92,3\% & 83,3\% & -          3 503 908  & -                97 263  & -         136,87  & -      2 354,34  & 97,22\% \\
          &       & RF    &       & 92,9\% & 83,3\% & -          2 850 198  &                           0   & -         120,28  &                     0    & 100,00\% \\
          &       & XGB   &       & 93,4\% & 83,3\% & -          3 151 393  &                           0    & -         124,60  &                     0    & 100,00\% \\
    \midrule
    \multirow{3}[2]{*}{A-16} & \multirow{3}[2]{*}{5096} & CART  & \multirow{3}[2]{*}{84,46\%} & 92,3\% & 83,8\% & -          3 778 933  & -             734 773  & -         147,49  & -         418,58  & 80,56\% \\
          &       & RF    &       & 92,9\% & 83,5\% & -          2 513 261  &                     8 914  & -         106,03  &              20,13  & 100,35\% \\
          &       & XGB   &       & 93,4\% & 83,6\% & -          2 920 405  &                  34 492  & -         115,47  &              45,75  & 101,18\% \\
    \bottomrule
    \end{tabular}%
\end{adjustbox}   

\begin{adjustbox}{width=1\textwidth}
    \begin{tabular}{cclcrrrrrrr}
    \multirow{2}[1]{*}{N°} & \multirow{2}[1]{*}{time (s)} & \multicolumn{1}{c}{\multirow{2}[1]{*}{Model}} & \multicolumn{1}{c}{\multirow{2}[1]{*}{\% target diff}} & \multicolumn{2}{c}{Accuracy} & \multicolumn{2}{c}{ Retention gain } & \multicolumn{2}{c}{RG/target} & \multicolumn{1}{c}{\multirow{2}[1]{*}{Improvement\footref{improvement_note}}} \\
          &       &       &       & \multicolumn{1}{c}{$y_i$} & \multicolumn{1}{c}{$\tilde{y}_i$} & \multicolumn{1}{c}{ $y_i$ } & \multicolumn{1}{c}{ $\tilde{y}_i$ } & \multicolumn{1}{c}{ $y_i$ } & \multicolumn{1}{c}{ $\tilde{y}_i$ } &  \\
    \midrule
    \multirow{3}[2]{*}{A-17} & \multirow{3}[2]{*}{5390} & CART  & \multirow{3}[2]{*}{28,74\%} & 92,3\% & 89,5\% &             5 100 456  &             4 899 479  &            199,11  &            279,88  & -3,94\% \\
          &       & RF    &       & 92,9\% & 89,8\% &             4 635 482  &             4 226 648  &            195,69  &            276,06  & -8,82\% \\
          &       & XGB   &       & 93,4\% & 90,2\% &             5 196 736  &             5 138 253  &            205,40  &            299,27  & -1,13\% \\
    \midrule
    \multirow{3}[2]{*}{A-18} & \multirow{3}[2]{*}{6452} & CART  & \multirow{3}[2]{*}{12,12\%} & 92,3\% & 91,3\% &          52 090 240  &          47 706 070  &        2 034,15  &        2 170,64  & -8,42\% \\
          &       & RF    &       & 92,9\% & 91,9\% &          46 171 160  &          42 049 900  &        1 949,05  &        2 082,36  & -8,93\% \\
          &       & XGB   &       & 93,4\% & 92,5\% &          51 629 950  &          52 606 740  &        2 040,95  &        2 339,70  & 1,89\% \\
    \midrule
    \multirow{3}[2]{*}{A-19} & \multirow{3}[2]{*}{4913} & CART  & \multirow{3}[2]{*}{64,89\%} & 92,3\% & 85,2\% &             2 798 173  &             3 182 143  &            109,11  &            481,60  & 13,72\% \\
          &       & RF    &       & 92,9\% & 85,2\% &             2 502 903  &             2 743 070  &            105,69  &            554,76  & 9,60\% \\
          &       & XGB   &       & 93,4\% & 85,6\% &             2 920 720  &             3 438 303  &            115,40  &            576,64  & 17,72\% \\
    \midrule
    \multirow{3}[2]{*}{A-20} & \multirow{3}[2]{*}{6160} & CART  & \multirow{3}[2]{*}{29,03\%} & 92,3\% & 89,6\% &          49 787 960  &          45 366 730  &        1 944,15  &        2 616,32  & -8,88\% \\
          &       & RF    &       & 92,9\% & 90,0\% &          44 038 580  &          39 947 830  &        1 859,05  &        2 547,89  & -9,29\% \\
          &       & XGB   &       & 93,4\% & 90,6\% &          49 353 940  &          49 789 670  &        1 950,95  &        2 796,17  & 0,88\% \\
    \midrule
    \multirow{3}[2]{*}{A-21} & \multirow{3}[2]{*}{5079} & CART  & \multirow{3}[2]{*}{51,69\%} & 92,3\% & 86,8\% &                482 682  &                544 887  &              18,85  &              53,99  & 12,89\% \\
          &       & RF    &       & 92,9\% & 86,8\% &                557 090  &                554 195  &              23,52  &              65,17  & -0,52\% \\
          &       & XGB   &       & 93,4\% & 87,1\% &                607 670  &                624 556  &              24,01  &              64,79  & 2,78\% \\
    \midrule
    \multirow{3}[2]{*}{A-22} & \multirow{3}[2]{*}{6199} & CART  & \multirow{3}[2]{*}{23,94\%} & 92,3\% & 90,2\% &             9 335 438  &             8 527 444  &            364,60  &            454,78  & -8,66\% \\
          &       & RF    &       & 92,9\% & 90,6\% &             8 570 307  &             7 931 029  &            361,80  &            460,42  & -7,46\% \\
          &       & XGB   &       & 93,4\% & 91,2\% &             9 518 466  &             9 581 934  &            376,27  &            501,56  & 0,67\% \\
    \midrule
    \multirow{3}[2]{*}{A-23} & \multirow{3}[2]{*}{4601} & CART  & \multirow{3}[2]{*}{89,51\%} & 92,3\% & 83,6\% & -          1 819 600  &                135 305  & -           71,15  &            121,80  & 107,44\% \\
          &       & RF    &       & 92,9\% & 83,5\% & -          1 575 489  &                159 620  & -           66,48  &            215,65  & 110,13\% \\
          &       & XGB   &       & 93,4\% & 83,7\% & -          1 668 346  &                228 226  & -           65,99  &            208,69  & 113,68\% \\
    \midrule
    \multirow{3}[2]{*}{A-24} & \multirow{3}[2]{*}{5650} & CART  & \multirow{3}[2]{*}{50,83\%} & 92,3\% & 87,0\% &             7 033 156  &             7 124 100  &            274,60  &            680,08  & 1,29\% \\
          &       & RF    &       & 92,9\% & 87,0\% &             6 437 729  &             6 364 477  &            271,80  &            711,89  & -1,14\% \\
          &       & XGB   &       & 93,4\% & 87,4\% &             7 242 450  &             7 840 770  &            286,27  &            771,71  & 8,26\% \\
    \midrule
    \multirow{3}[2]{*}{A-25} & \multirow{3}[2]{*}{5379} & CART  & \multirow{3}[2]{*}{30,97\%} & 92,3\% & 89,2\% &             4 160 423  &             3 882 623  &            162,44  &            241,06  & -6,68\% \\
          &       & RF    &       & 92,9\% & 89,5\% &             4 018 432  &             3 666 219  &            169,65  &            249,54  & -8,76\% \\
          &       & XGB   &       & 93,4\% & 90,0\% &             4 455 108  &             4 410 629  &            176,09  &            267,87  & -1,00\% \\
    \midrule
    \multirow{3}[2]{*}{A-26} & \multirow{3}[2]{*}{6410} & CART  & \multirow{3}[2]{*}{12,52\%} & 92,3\% & 91,3\% &          49 612 660  &          45 948 690  &        1 937,51  &        2 083,30  & -7,39\% \\
          &       & RF    &       & 92,9\% & 91,9\% &          44 548 720  &          40 814 960  &        1 880,59  &        2 029,68  & -8,38\% \\
          &       & XGB   &       & 93,4\% & 92,5\% &          49 676 000  &          50 549 740  &        1 963,72  &        2 260,20  & 1,76\% \\
    \midrule
    \multirow{3}[2]{*}{A-27} & \multirow{3}[2]{*}{4887} & CART  & \multirow{3}[2]{*}{66,67\%} & 92,3\% & 85,1\% &             1 858 140  &             2 575 538  &              72,44  &            442,86  & 38,61\% \\
          &       & RF    &       & 92,9\% & 85,0\% &             1 885 853  &             2 387 018  &              79,65  &            531,25  & 26,57\% \\
          &       & XGB   &       & 93,4\% & 85,4\% &             2 179 093  &             2 879 880  &              86,09  &            544,35  & 32,16\% \\
    \midrule
    \multirow{3}[2]{*}{A-28} & \multirow{3}[2]{*}{6047} & CART  & \multirow{3}[2]{*}{29,42\%} & 92,3\% & 89,4\% &          47 310 370  &          43 168 880  &        1 847,51  &        2 519,41  & -8,75\% \\
          &       & RF    &       & 92,9\% & 89,9\% &          42 416 140  &          38 573 620  &        1 790,59  &        2 504,61  & -9,06\% \\
          &       & XGB   &       & 93,4\% & 90,5\% &          47 399 990  &          47 812 830  &        1 873,72  &        2 721,63  & 0,87\% \\
    \midrule
    \multirow{3}[2]{*}{A-29} & \multirow{3}[2]{*}{5070} & CART  & \multirow{3}[2]{*}{53,79\%} & 92,3\% & 86,5\% & -             204 467  & -                  5 098  & -              7,95  & -              1,66  & 97,51\% \\
          &       & RF    &       & 92,9\% & 86,1\% &                163 014  &                273 435  &                6,90  &              40,30  & 67,74\% \\
          &       & XGB   &       & 93,4\% & 86,8\% &                115 297  &                248 982  &                4,55  &              28,64  & 115,95\% \\
    \midrule
    \multirow{3}[2]{*}{A-30} & \multirow{3}[2]{*}{6179} & CART  & \multirow{3}[2]{*}{24,36\%} & 92,3\% & 90,3\% &             7 522 978  &             7 058 487  &            293,94  &            382,06  & -6,17\% \\
          &       & RF    &       & 92,9\% & 90,6\% &             7 534 275  &             7 068 293  &            318,08  &            411,80  & -6,18\% \\
          &       & XGB   &       & 93,4\% & 91,2\% &             8 219 857  &             8 265 167  &            324,94  &            442,88  & 0,55\% \\
    \midrule
    \multirow{3}[2]{*}{A-31} & \multirow{3}[2]{*}{4627} & CART  & \multirow{3}[2]{*}{90,18\%} & 92,3\% & 83,6\% & -          2 506 749  & -             139 983  & -           97,95  & -         121,44  & 94,42\% \\
          &       & RF    &       & 92,9\% & 83,5\% & -          1 969 564  &                  73 101  & -           83,10  &            111,49  & 103,71\% \\
          &       & XGB   &       & 93,4\% & 83,6\% & -          2 160 719  &                  76 641  & -           85,45  &              93,28  & 103,55\% \\
    \midrule
    \multirow{3}[2]{*}{A-32} & \multirow{3}[2]{*}{5679} & CART  & \multirow{3}[2]{*}{51,25\%} & 92,3\% & 86,8\% &             5 220 695  &             5 811 833  &            203,94  &            583,55  & 11,32\% \\
          &       & RF    &       & 92,9\% & 86,9\% &             5 401 696  &             5 269 505  &            228,08  &            605,69  & -2,45\% \\
          &       & XGB   &       & 93,4\% & 87,4\% &             5 943 841  &             6 682 230  &            234,94  &            670,03  & 12,42\% \\
    \bottomrule
    \end{tabular}%
 \end{adjustbox}       
\newpage
\begin{adjustbox}{width=1\textwidth}
    \begin{tabular}{cclcrrrrrrr}
    \multirow{2}[1]{*}{N°} & \multirow{2}[1]{*}{time (s)} & \multicolumn{1}{c}{\multirow{2}[1]{*}{Model}} & \multicolumn{1}{c}{\multirow{2}[1]{*}{\% target diff}} & \multicolumn{2}{c}{Accuracy} & \multicolumn{2}{c}{ Retention gain } & \multicolumn{2}{c}{ RG/target } & \multicolumn{1}{c}{\multirow{2}[1]{*}{Improvement\footref{improvement_note}}} \\
          &       &       &       & \multicolumn{1}{c}{$y_i$} & \multicolumn{1}{c}{$\tilde{y}_i$} & \multicolumn{1}{c}{ $y_i$ } & \multicolumn{1}{c}{ $\tilde{y}_i$ } & \multicolumn{1}{c}{ $y_i$ } & \multicolumn{1}{c}{ $\tilde{y}_i$ } &  \\
    \midrule
    \multirow{3}[2]{*}{B-1} & \multirow{3}[2]{*}{4778} & CART  & \multirow{3}[2]{*}{75,89\%} & 92,3\% & 84,0\% & -             627 165  & -             148 913  & -           24,46  & -           65,19  & 76,26\% \\
          &       & RF    &       & 92,9\% & 83,7\% & -             280 855  &                  11 973  & -           11,84  &                9,57  & 104,26\% \\
          &       & XGB   &       & 93,4\% & 84,1\% & -             366 103  &                  25 099  & -           14,47  &              12,30  & 106,86\% \\
    \midrule
    \multirow{3}[2]{*}{B-2} & \multirow{3}[2]{*}{6074} & CART  & \multirow{3}[2]{*}{29,70\%} & 92,3\% & 89,7\% &             3 862 156  &             3 397 247  &            150,95  &            203,11  & -12,04\% \\
          &       & RF    &       & 92,9\% & 89,9\% &             4 127 224  &             3 550 730  &            174,26  &            230,67  & -13,97\% \\
          &       & XGB   &       & 93,4\% & 90,6\% &             4 451 686  &             4 408 819  &            175,99  &            250,17  & -0,96\% \\
    \midrule
    \multirow{3}[2]{*}{B-3} & \multirow{3}[2]{*}{4528} & CART  & \multirow{3}[2]{*}{96,60\%} & 92,3\% & 83,2\% & -          2 929 448  & -                85 465  & -         114,46  & -      1 482,06  & 97,08\% \\
          &       & RF    &       & 92,9\% & 83,3\% & -          2 413 433  &                     3 724  & -         101,84  & -         108,33  & 100,15\% \\
          &       & XGB   &       & 93,4\% & 83,3\% & -          2 642 119  &                     9 092  & -         104,47  &              93,79  & 100,34\% \\
    \midrule
    \multirow{3}[2]{*}{B-4} & \multirow{3}[2]{*}{5476} & CART  & \multirow{3}[2]{*}{60,93\%} & 92,3\% & 85,9\% &             1 559 874  &             2 471 262  &              60,95  &            329,63  & 58,43\% \\
          &       & RF    &       & 92,9\% & 85,8\% &             1 994 645  &             2 517 111  &              84,26  &            422,45  & 26,19\% \\
          &       & XGB   &       & 93,4\% & 86,3\% &             2 175 670  &             3 089 897  &              85,99  &            422,77  & 42,02\% \\
    \midrule
    \multirow{3}[2]{*}{B-5} & \multirow{3}[2]{*}{4708} & CART  & \multirow{3}[2]{*}{84,33\%} & 92,3\% & 83,4\% & -             857 439  & -             159 856  & -           33,45  & -         218,16  & 81,36\% \\
          &       & RF    &       & 92,9\% & 83,3\% & -             484 459  &                          40  & -           20,44  &                7,23  & 100,01\% \\
          &       & XGB   &       & 93,4\% & 83,3\% & -             596 203  &                        897  & -           23,57  &              46,96  & 100,15\% \\
    \midrule
    \multirow{3}[2]{*}{B-6} & \multirow{3}[2]{*}{5906} & CART  & \multirow{3}[2]{*}{36,63\%} & 92,3\% & 88,8\% &                705 721  &                922 490  &              27,69  &              60,21  & 30,72\% \\
          &       & RF    &       & 92,9\% & 88,9\% &             1 352 182  &             1 269 349  &              57,11  &              97,63  & -6,13\% \\
          &       & XGB   &       & 93,4\% & 89,6\% &             1 342 882  &             1 428 722  &              53,09  &              96,76  & 6,39\% \\
    \midrule
    \multirow{3}[2]{*}{B-7} & \multirow{3}[2]{*}{4400} & CART  & \multirow{3}[2]{*}{98,49\%} & 92,3\% & 83,2\% & -          3 159 722  & -                39 633  & -         123,45  & -      1 230,61  & 98,75\% \\
          &       & RF    &       & 92,9\% & 83,3\% & -          2 617 037  &                     1 024  & -         110,44  &                0,56  & 100,04\% \\
          &       & XGB   &       & 93,4\% & 83,3\% & -          2 872 219  &                        295  & -         113,57  &              19,31  & 100,01\% \\
    \midrule
    \multirow{3}[2]{*}{B-8} & \multirow{3}[2]{*}{5278} & CART  & \multirow{3}[2]{*}{73,18\%} & 92,3\% & 84,6\% & -          1 596 562  &                169 852  & -           62,31  &              41,78  & 110,64\% \\
          &       & RF    &       & 92,9\% & 84,6\% & -             780 396  &                637 625  & -           32,89  &            194,52  & 181,71\% \\
          &       & XGB   &       & 93,4\% & 85,0\% & -             933 133  &                780 845  & -           36,91  &            188,79  & 183,68\% \\
    \midrule
    \multirow{3}[2]{*}{B-9} & \multirow{3}[2]{*}{4601} & CART  & \multirow{3}[2]{*}{94,12\%} & 92,3\% & 83,3\% & -          2 380 789  & -             113 444  & -           92,86  & -         840,25  & 95,24\% \\
          &       & RF    &       & 92,9\% & 83,3\% & -          1 403 468  &                        317  & -           59,21  &                7,96  & 100,02\% \\
          &       & XGB   &       & 93,4\% & 83,3\% & -          1 724 731  &                     3 980  & -           68,17  &            149,44  & 100,23\% \\
    \midrule
    \multirow{3}[2]{*}{B-10} & \multirow{3}[2]{*}{5947} & CART  & \multirow{3}[2]{*}{35,98\%} & 92,3\% & 89,0\% & -             760 449  &                429 196  & -           29,35  &              29,80  & 156,44\% \\
          &       & RF    &       & 92,9\% & 88,5\% &             1 175 540  &             1 354 131  &              49,71  &            118,11  & 15,19\% \\
          &       & XGB   &       & 93,4\% & 89,8\% &                871 455  &             1 456 080  &              34,48  &              96,25  & 67,09\% \\
    \midrule
    \multirow{3}[2]{*}{B-11} & \multirow{3}[2]{*}{4229} & CART  & \multirow{3}[2]{*}{99,16\%} & 92,3\% & 83,3\% & -          4 683 072  & -                48 985  & -         182,86  & -      1 186,22  & 98,95\% \\
          &       & RF    &       & 92,9\% & 83,3\% & -          3 536 046  &                           0    & -         149,21  &                     0    & 100,00\% \\
          &       & XGB   &       & 93,4\% & 83,3\% & -          4 000 747  &                           0    & -         158,17  &                     0    & 100,00\% \\
    \midrule
    \multirow{3}[2]{*}{B-12} & \multirow{3}[2]{*}{5391} & CART  & \multirow{3}[2]{*}{66,76\%} & 92,3\% & 85,0\% & -          3 062 732  & -             388 289  & -         119,35  & -           80,44  & 87,32\% \\
          &       & RF    &       & 92,9\% & 84,7\% & -             957 039  &                710 688  & -           40,29  &            220,55  & 174,26\% \\
          &       & XGB   &       & 93,4\% & 85,3\% & -          1 404 561  &                834 198  & -           55,52  &            163,88  & 159,39\% \\
    \midrule
    \multirow{3}[2]{*}{B-13} & \multirow{3}[2]{*}{4493} & CART  & \multirow{3}[2]{*}{96,30\%} & 92,3\% & 83,3\% & -          2 358 179  & -             159 922  & -           91,98  & -      2 793,13  & 93,22\% \\
          &       & RF    &       & 92,9\% & 83,3\% & -          1 384 098  &                           0    & -           58,40  &                     0    & 100,00\% \\
          &       & XGB   &       & 93,4\% & 83,3\% & -          1 705 577  &                           0    & -           67,42  &                     0    & 100,00\% \\
    \midrule
    \multirow{3}[2]{*}{B-14} & \multirow{3}[2]{*}{5851} & CART  & \multirow{3}[2]{*}{42,98\%} & 92,3\% & 87,8\% & -          3 251 762  & -          1 761 821  & -         126,63  & -         143,20  & 45,82\% \\
          &       & RF    &       & 92,9\% & 86,4\% & -          1 013 089  &                  79 273  & -           42,69  &              11,90  & 107,82\% \\
          &       & XGB   &       & 93,4\% & 83,3\% & -          1 582 006  &                     4 396  & -           62,52  &            287,68  & 100,28\% \\
    \midrule
    \multirow{3}[2]{*}{B-15} & \multirow{3}[2]{*}{4040} & CART  & \multirow{3}[2]{*}{99,67\%} & 92,3\% & 83,3\% & -          4 660 462  & -                38 969  & -         181,98  & -      2 075,03  & 99,16\% \\
          &       & RF    &       & 92,9\% & 83,3\% & -          3 516 676  &                           0    & -         148,40  &                     0    & 100,00\% \\
          &       & XGB   &       & 93,4\% & 83,3\% & -          3 981 592  &                        161  & -         157,42  &              10,53  & 100,00\% \\
    \midrule
    \multirow{3}[2]{*}{B-16} & \multirow{3}[2]{*}{5182} & CART  & \multirow{3}[2]{*}{77,97\%} & 92,3\% & 84,2\% & -          5 554 044  & -          1 491 522  & -         216,63  & -         549,23  & 73,15\% \\
          &       & RF    &       & 92,9\% & 83,6\% & -          3 145 668  &                  52 475  & -         132,69  &              84,54  & 101,67\% \\
          &       & XGB   &       & 93,4\% & 83,3\% & -          3 858 022  &                           0    & -         152,52  &                     0    & 100,00\% \\
    \bottomrule
    \end{tabular}%
\end{adjustbox}        
\newpage
\begin{adjustbox}{width=1\textwidth}
    \begin{tabular}{cclcrrrrrrr}
    \multirow{2}[1]{*}{N°} & \multirow{2}[1]{*}{time (s)} & \multicolumn{1}{c}{\multirow{2}[1]{*}{Model}} & \multicolumn{1}{c}{\multirow{2}[1]{*}{\% target diff}} & \multicolumn{2}{c}{Accuracy} & \multicolumn{2}{c}{ Retention gain } & \multicolumn{2}{c}{RG/target} & \multicolumn{1}{c}{\multirow{2}[1]{*}{Improvement\footref{improvement_note}}} \\
          &       &       &       & \multicolumn{1}{c}{$y_i$} & \multicolumn{1}{c}{$\tilde{y}_i$} & \multicolumn{1}{c}{ $y_i$ } & \multicolumn{1}{c}{ $\tilde{y}_i$ } & \multicolumn{1}{c}{ $y_i$ } & \multicolumn{1}{c}{ $\tilde{y}_i$ } &  \\
    \midrule
    \multirow{3}[2]{*}{B-17} & \multirow{3}[2]{*}{5324} & CART  & \multirow{3}[2]{*}{32,66\%} & 92,3\% & 88,9\% &             3 361 471  &             3 037 200  &            131,25  &            191,31  & -9,65\% \\
          &       & RF    &       & 92,9\% & 89,3\% &             3 241 680  &             2 911 023  &            136,86  &            204,43  & -10,20\% \\
          &       & XGB   &       & 93,4\% & 89,6\% &             3 596 593  &             3 546 671  &            142,15  &            222,04  & -1,39\% \\
    \midrule
    \multirow{3}[2]{*}{B-18} & \multirow{3}[2]{*}{6411} & CART  & \multirow{3}[2]{*}{13,83\%} & 92,3\% & 91,1\% &          39 860 670  &          37 695 680  &        1 556,66  &        1 778,71  & -5,43\% \\
          &       & RF    &       & 92,9\% & 91,7\% &          35 787 050  &          32 345 100  &        1 510,72  &        1 654,32  & -9,62\% \\
          &       & XGB   &       & 93,4\% & 92,0\% &          39 908 670  &          40 886 810  &        1 577,61  &        1 848,71  & 2,45\% \\
    \midrule
    \multirow{3}[2]{*}{B-19} & \multirow{3}[2]{*}{4853} & CART  & \multirow{3}[2]{*}{70,34\%} & 92,3\% & 84,7\% &             1 059 189  &             1 813 631  &              41,25  &            392,14  & 71,23\% \\
          &       & RF    &       & 92,9\% & 84,8\% &             1 109 101  &             1 808 616  &              46,86  &            474,33  & 63,07\% \\
          &       & XGB   &       & 93,4\% & 85,0\% &             1 320 578  &             2 141 271  &              52,15  &            482,34  & 62,15\% \\
    \midrule
    \multirow{3}[2]{*}{B-20} & \multirow{3}[2]{*}{5973} & CART  & \multirow{3}[2]{*}{31,76\%} & 92,3\% & 89,2\% &          37 558 390  &          34 068 550  &        1 466,66  &        2 125,97  & -9,29\% \\
          &       & RF    &       & 92,9\% & 89,4\% &          33 654 470  &          30 032 580  &        1 420,72  &        2 072,47  & -10,76\% \\
          &       & XGB   &       & 93,4\% & 90,1\% &          37 632 650  &          38 008 480  &        1 487,61  &        2 277,17  & 1,00\% \\
    \midrule
    \multirow{3}[2]{*}{B-21} & \multirow{3}[2]{*}{5228} & CART  & \multirow{3}[2]{*}{41,79\%} & 92,3\% & 87,7\% &             1 136 879  &             1 179 837  &              44,40  &              92,50  & 3,78\% \\
          &       & RF    &       & 92,9\% & 88,1\% &             1 276 808  &             1 188 256  &              53,91  &            104,81  & -6,94\% \\
          &       & XGB   &       & 93,4\% & 88,7\% &             1 385 145  &             1 356 864  &              54,74  &            104,76  & -2,04\% \\
    \midrule
    \multirow{3}[2]{*}{B-22} & \multirow{3}[2]{*}{6296} & CART  & \multirow{3}[2]{*}{19,52\%} & 92,3\% & 90,7\% &          18 704 980  &          17 177 190  &            730,55  &            852,81  & -8,17\% \\
          &       & RF    &       & 92,9\% & 91,1\% &          17 182 100  &          15 732 340  &            725,34  &            859,29  & -8,44\% \\
          &       & XGB   &       & 93,4\% & 91,5\% &          19 071 370  &          19 050 020  &            753,90  &            939,00  & -0,11\% \\
    \midrule
    \multirow{3}[2]{*}{B-23} & \multirow{3}[2]{*}{4746} & CART  & \multirow{3}[2]{*}{81,36\%} & 92,3\% & 84,1\% & -          1 165 404  &                458 223  & -           45,60  &            172,83  & 139,32\% \\
          &       & RF    &       & 92,9\% & 84,0\% & -             855 770  &                525 335  & -           36,09  &            288,55  & 161,39\% \\
          &       & XGB   &       & 93,4\% & 84,1\% & -             890 871  &                645 445  & -           35,26  &            310,86  & 172,45\% \\
    \midrule
    \multirow{3}[2]{*}{B-24} & \multirow{3}[2]{*}{5845} & CART  & \multirow{3}[2]{*}{40,47\%} & 92,3\% & 88,2\% &          16 402 700  &          15 013 310  &            640,55  &        1 093,43  & -8,47\% \\
          &       & RF    &       & 92,9\% & 88,4\% &          15 049 520  &          13 423 040  &            635,34  &        1 122,81  & -10,81\% \\
          &       & XGB   &       & 93,4\% & 88,9\% &          16 795 360  &          17 144 260  &            663,90  &        1 247,50  & 2,08\% \\
    \midrule
    \multirow{3}[2]{*}{B-25} & \multirow{3}[2]{*}{5274} & CART  & \multirow{3}[2]{*}{37,42\%} & 92,3\% & 88,6\% &             1 607 847  &             1 839 864  &              62,84  &            126,33  & 14,43\% \\
          &       & RF    &       & 92,9\% & 88,7\% &             2 119 067  &             1 923 982  &              89,49  &            152,71  & -9,21\% \\
          &       & XGB   &       & 93,4\% & 89,2\% &             2 237 965  &             2 194 469  &              88,45  &            155,54  & -1,94\% \\
    \midrule
    \multirow{3}[2]{*}{B-26} & \multirow{3}[2]{*}{6425} & CART  & \multirow{3}[2]{*}{14,83\%} & 92,3\% & 91,1\% &          35 238 060  &          32 690 970  &        1 376,37  &        1 558,26  & -7,23\% \\
          &       & RF    &       & 92,9\% & 91,6\% &          32 835 370  &          29 986 540  &        1 386,17  &        1 543,12  & -8,68\% \\
          &       & XGB   &       & 93,4\% & 92,0\% &          36 328 440  &          36 803 630  &        1 436,10  &        1 688,53  & 1,31\% \\
    \midrule
    \multirow{3}[2]{*}{B-27} & \multirow{3}[2]{*}{4811} & CART  & \multirow{3}[2]{*}{73,92\%} & 92,3\% & 84,3\% & -             694 436  &                751 404  & -           27,16  &            226,99  & 208,20\% \\
          &       & RF    &       & 92,9\% & 84,4\% & -                13 512  &             1 018 369  & -              0,51  &            356,48  & 7636,98\% \\
          &       & XGB   &       & 93,4\% & 84,7\% & -                38 050  &             1 253 252  & -              1,55  &            345,94  & 3393,68\% \\
    \midrule
    \multirow{3}[2]{*}{B-28} & \multirow{3}[2]{*}{5995} & CART  & \multirow{3}[2]{*}{32,61\%} & 92,3\% & 89,1\% &          32 935 780  &          29 342 930  &        1 286,37  &        1 847,71  & -10,91\% \\
          &       & RF    &       & 92,9\% & 89,4\% &          30 702 790  &          27 725 620  &        1 296,17  &        1 933,38  & -9,70\% \\
          &       & XGB   &       & 93,4\% & 90,0\% &          34 052 420  &          34 390 060  &        1 346,10  &        2 094,90  & 0,99\% \\
    \midrule
    \multirow{3}[2]{*}{B-29} & \multirow{3}[2]{*}{5143} & CART  & \multirow{3}[2]{*}{47,03\%} & 92,3\% & 87,3\% & -             363 861  &                  55 985  & -           14,12  &                3,38  & 115,39\% \\
          &       & RF    &       & 92,9\% & 87,4\% &                377 170  &                488 284  &              15,95  &              49,62  & 29,46\% \\
          &       & XGB   &       & 93,4\% & 88,0\% &                275 772  &                491 567  &              10,89  &              44,89  & 78,25\% \\
    \midrule
    \multirow{3}[2]{*}{B-30} & \multirow{3}[2]{*}{6243} & CART  & \multirow{3}[2]{*}{20,47\%} & 92,3\% & 90,7\% &          14 747 500  &          13 838 380  &            576,23  &            690,22  & -6,16\% \\
          &       & RF    &       & 92,9\% & 91,1\% &          14 816 830  &          13 378 460  &            625,54  &            743,34  & -9,71\% \\
          &       & XGB   &       & 93,4\% & 91,5\% &          16 146 490  &          16 169 440  &            638,30  &            814,80  & 0,14\% \\
    \midrule
    \multirow{3}[2]{*}{B-31} & \multirow{3}[2]{*}{4730} & CART  & \multirow{3}[2]{*}{83,83\%} & 92,3\% & 83,7\% & -          2 666 144  & -             487 716  & -         104,12  & -         267,75  & 81,71\% \\
          &       & RF    &       & 92,9\% & 83,7\% & -          1 755 409  &                139 545  & -           74,05  &            102,66  & 107,95\% \\
          &       & XGB   &       & 93,4\% & 83,7\% & -          2 000 244  &                134 199  & -           79,11  &            130,13  & 106,71\% \\
    \midrule
    \multirow{3}[2]{*}{B-32} & \multirow{3}[2]{*}{5865} & CART  & \multirow{3}[2]{*}{41,41\%} & 92,3\% & 88,0\% &          12 445 210  &          11 693 070  &            486,23  &            884,49  & -6,04\% \\
          &       & RF    &       & 92,9\% & 88,3\% &          12 684 250  &          11 381 260  &            535,54  &            971,28  & -10,27\% \\
          &       & XGB   &       & 93,4\% & 88,8\% &          13 870 470  &          14 101 470  &            548,30  &        1 048,38  & 1,67\% \\
    \bottomrule
    \end{tabular}%
\end{adjustbox}    
\footnotetext{In order to account for negative retention gains, the improvement is computed with an absolute value for the denominator. This leads to a rather unintuitive improvement measure whenever one of the models yields negative RG and the other positive RG.\label{improvement_note}}

\end{document}